%% file: Cai-Guo-Confidence_interval_AoS_revised_version__Main_paper_11-25-2015.tex
\documentclass[aos,preprint]{imsart}

\RequirePackage{amsthm,amsmath}
\RequirePackage[numbers]{natbib}%
\RequirePackage[colorlinks,citecolor=blue,urlcolor=blue]{hyperref}

\usepackage{amsgen,amsmath,amstext,amsbsy,amsopn,amssymb,latexsym,comment}
\usepackage{array}
\usepackage{cite}
\usepackage{multirow}
\usepackage{graphicx}
\usepackage{amssymb}
\usepackage{bibunits}
\usepackage[dvipsnames,table,dvipsnames*, svgnames*, hyperref]{xcolor}
\arxiv{arXiv:0000.0000}

\newcommand{\red }{\color{red}}

\newcommand{\beq}{\begin{equation}}
\newcommand{\eeq}{\end{equation}}
\newcommand{\beas}{\begin{eqnarray*}}
\newcommand{\eeas}{\end{eqnarray*}}
\newcommand{\bea}{\begin{eqnarray}}
\newcommand{\eea}{\end{eqnarray}}
\newcommand{\bei}{\begin{itemize}}
\newcommand{\eei}{\end{itemize}}
\newcommand{\ben}{\begin{enumerate}}
\newcommand{\een}{\end{enumerate}}
\newcommand{\argmin}{\mathop{\rm arg\min}}

\newtheorem{Corollary}{Corollary}
\newtheorem{Proposition}{Proposition}
\newtheorem{Lemma}{Lemma}

\newtheorem{Theorem}{Theorem}

\newtheorem{Remark}{Remark}

\newcommand{\za}{z_{\alpha/2}}
\newcommand{\zas}{z_{\alpha_0/2}}

\newcommand{\CIZa}{{\rm CI}_{\alpha}\left(\beta_1,Z\right)}
\newcommand{\CIZb}{{\rm CI}_{\alpha}\left(\sum \beta_i,Z\right)}
\newcommand{\CIZl}{{\rm CI}_{\alpha}\left(\xi^{\intercal} \beta, Z\right)}
\newcommand{\CIZs}{{\rm CI}^{S}_{\alpha}\left(\xi^{\intercal} \beta, Z\right)}
\newcommand{\CIZd}{{\rm CI}^{D}_{\alpha}\left(\xi^{\intercal} \beta, Z\right)}
\newcommand{\CIZk}{{\rm CI}^{{\rm I}}_{\alpha}\left(\xi^{\intercal} \beta, Z\right)}

\newcommand{\CIZg}{{\rm CI}_{\alpha}\left(\Tb,Z\right)}
\newcommand{\CIZtrans}{{\rm CI}_{\alpha}\left(\psi_1,Z\right)}
\newcommand{\SCIZ}{ \II_{\alpha} \left(\Theta,\Tb\right)}
\newcommand{\SCIZa}{\II_{\alpha} \left(\Theta\left(k\right),\beta_1\right)}
\newcommand{\SCIZb}{\II_{\alpha} \left(\Theta\left(k\right),\sum\beta_i\right)}
\newcommand{\SCIZl}{\II_{\alpha} \left(\Theta\left(k\right),\xi^{\intercal}\beta\right)}
\newcommand{\SCIZg}{\II_{\alpha} \left(\Theta\left(k\right),\Tb\right)}
\newcommand{\SCIZtrans}{\II_{\alpha} \left(\GG\left(k\right),\psi_1\right)}
\newcommand{\suppx}{{\rm supp}(x)}
\newcommand{\suppxi}{{\rm supp}(\xi)}

\newcommand{\esparse}{k\ll \frac{\sqrt{n}}{\log p}}
\newcommand{\usparse}{k\lesssim \frac{\sqrt{n}}{\log p}}
\newcommand{\msparse}{\frac{\sqrt{n}}{\log p}\ll k\lesssim \frac{n}{\log p}}

\newcommand{\Tf}{{\rm T}\left(\beta\right)}
\newcommand{\Tb}{{\rm T}}
\newcommand{\linfun}{\xi^{\intercal} \beta}
\newcommand{\liminfnp}{\liminf\limits_{n,p \rightarrow \infty}}
\newcommand{\Xj}{X_{\cdot j}}

\newcommand{\Md}{M_1}
\newcommand{\Mv}{M_2}

\newcommand{\TV}{{\rm TV}}

\newcommand{\alphab}{\boldsymbol{\delta}}

\newcommand{\xa}{X^{\left(1\right)}}
\newcommand{\xb}{X^{\left(2\right)}}

\newcommand{\ya}{y^{\left(1\right)}}
\newcommand{\yb}{y^{\left(2\right)}}

\newcommand{\uu}{u}
\newcommand{\R}{\mathbb{R}}
\newcommand{\E}{{\mathbb{E}}}
\newcommand{\PP}{{\mathbb{P}_{\theta}}}
\newcommand{\FF}{{\mathcal{F}}}
\newcommand{\GG}{{\mathcal{G}}}
\newcommand{\HH}{{\mathcal{H}}}
\newcommand{\II}{{\mathcal{I}}}

\newcommand{\T}{{\rm T}}
\newcommand{\supp}{{\rm supp}}

\startlocaldefs
\numberwithin{equation}{section}
\theoremstyle{plain}

\endlocaldefs

\begin{document}

\begin{frontmatter}
\title{Confidence Intervals for High-Dimensional Linear Regression: Minimax Rates and Adaptivity\thanksref{T1}}
\runtitle{High-dimensional Confidence Interval}
\thankstext{T1}{The research was supported in part by NSF Grants DMS-1208982 and DMS-1403708, and NIH Grant R01 CA127334..}

\begin{aug}
\author{\fnms{T. Tony} \snm{Cai}\ead[label=e1]{tcai@wharton.upenn.edu}},
\and
\author{\fnms{Zijian} \snm{Guo}\ead[label=e2]{zijguo@wharton.upenn.edu}
\ead[label=u1,url]{URL: http://www-stat.wharton.upenn.edu/$\sim$tcai/}}
\runauthor{T. T. Cai and Z. Guo}
\affiliation{University of Pennsylvania}
\address{DEPARTMENT OF STATISTICS\\
THE WHARTON SCHOOL\\
UNIVERSITY OF PENNSYLVANIA\\
PHILADELPHIA, PENNSYLVANIA 19104\\
USA\\
\printead{e1}\\
\phantom{E-mail:\ }\printead*{e2}\\
\printead*{u1}\phantom{URL:\ }\\
}

\end{aug}

\begin{abstract}
Confidence sets play a fundamental role in statistical inference. In this paper, we consider  confidence intervals for high dimensional linear regression with random design. We first establish the convergence rates of the minimax expected length for confidence intervals in the oracle setting where the sparsity parameter is given. The focus is then on the problem of adaptation to sparsity for the construction of confidence intervals. Ideally, an adaptive confidence interval should have its length automatically adjusted to the sparsity of the unknown regression vector, while maintaining a prespecified coverage probability.  It is shown that such a goal is in general not attainable, except when the sparsity parameter is restricted to a small region over which the confidence intervals have the optimal length of the usual parametric rate. It is further demonstrated that the lack of adaptivity is not due to the conservativeness of the minimax framework, but is fundamentally caused by the difficulty of learning the bias accurately.
\end{abstract}

\begin{keyword}[class=MSC]
\kwd[Primary ]{62G15}
\kwd[; secondary ]{62C20}\kwd{62H35}
\end{keyword}

\begin{keyword}
\kwd{Adaptivity, confidence interval, coverage probability, expected length, high-dimensional linear regression, minimaxity, sparsity.}
\end{keyword}
\end{frontmatter}
\input{CIMaintext}
\section*{Acknowledgements}
The authors thank Zhao Ren for the discussion on the confidence intervals for linear functionals with sparse loadings. \begin{supplement}
\stitle{Supplement to ``Confidence Intervals for High-Dimensional Linear Regression: Minimax Rates and Adaptivity''.}
\slink[url]{http://www-stat.wharton.upenn.edu/$\sim$tcai/paper/CI-Reg-Supplement.pdf}
\sdescription{
Detailed proofs of the adaptivity lower bound and minimax upper bound for confidence intervals of the linear functional $\xi^{\intercal} \beta$ with a dense loading $\xi$ are given. The minimax rates and adaptivity of confidence intervals of the linear functional $\xi^{\intercal} \beta$ are established when there is prior knowledge that $\Omega={\rm I}$ and $\sigma=\sigma_0$.  Extra propositions and technical lemmas are also proved in the supplement. 
}
\end{supplement}
\bibliographystyle{plain}
\bibliography{HighDimRef}

\end{document}

%% file: CIMaintext.tex

\section{Introduction}


Driven by a wide range of applications,  high-dimensional linear regression, where the dimension $p$ can be much larger than the sample size $n$, has received significant recent attention. The linear model is
\begin{equation}
y=X\beta+\epsilon, \quad \epsilon \sim N(0,\sigma^2 {\rm I}),
\label{eq: original linear model}
\end{equation}
where $y\in \R^{n}$, $X\in \R^{n\times p}$, and $\beta\in \R^p$. Several penalized/constrained $\ell_1$ minimization methods, including the Lasso \citep{tibshirani1996regression}, Dantzig Selector \citep{candes2007dantzig}, square-root Lasso \citep{belloni2011square},
and scaled Lasso \citep{sun2012scaled} have been proposed and studied. Under regularity conditions on the design matrix $X$,  these methods with a suitable choice of the tuning parameter have been shown to achieve the optimal rate of convergence $k \frac{\log p}{n}$ under the squared error loss over the set of $k$-sparse regression coefficient vectors with $k\leq c \frac{n}{\log p}$ where $c>0$ is a constant. That is, there exists some constant $C>0$ such that
\begin{equation}
\sup_{\|\beta\|_0\leq k} \mathbb{P}\left(\|\widehat{\beta}-\beta\|_2^2 > C k \frac{\log p}{n}\right)=o(1),
\label{eq: minimax estimation}
\end{equation}
where 
$\|\beta\|_0$ denotes the number of the nonzero coordinates of a vector $\beta\in \R^p$.
See, for example, \citep{verzelen2012minimax,bickel2009simultaneous,candes2007dantzig,sun2012scaled}. A key feature of the estimation problem is that the optimal rate can be achieved adaptively with respect to the sparsity parameter $k$.

Confidence sets play a fundamental role in statistical inference and confidence intervals for high-dimensional linear regression have been actively studied recently with a focus on inference for individual coordinates. But, compared to point estimation, there is still a paucity of methods and fundamental theoretical results on confidence intervals for high-dimensional regression. Zhang and Zhang \citep{zhang2014confidence} was the first to introduce the idea of de-biasing for constructing a valid confidence interval for a single coordinate $\beta_i$. The confidence interval is centered at a low-dimensional projection estimator obtained through bias correction via score vector using the scaled Lasso as the initial estimator.
\citep{javanmard2013hypothesis,javanmard2014confidence,van2014asymptotically} also used de-biasing for  the construction of confidence intervals and  \citep{van2014asymptotically} established asymptotic efficiency for the proposed estimator.
 All the aforementioned papers \citep{zhang2014confidence, javanmard2013hypothesis, javanmard2014confidence, van2014asymptotically} have focused on the ultra-sparse case where the sparsity $k \ll \frac{\sqrt{n}}{\log p}$ is assumed.
 Under such a sparsity condition, the expected length of the confidence intervals constructed in \citep{zhang2014confidence, javanmard2014confidence, van2014asymptotically}  is at the parametric rate $\frac{1}{\sqrt{n}}$ and the procedures do not depend on the specific value of $k$.

Compared to point estimation where the sparsity condition $k \ll \frac{n}{\log p}$ is sufficient for estimation consistency (see equation \eqref{eq: minimax estimation}), the condition $\esparse$ for valid confidence intervals is much stronger. There are several natural questions: What happens in the region where  ${\sqrt{n}\over \log p} \lesssim k \lesssim {n \over \log p}$?
Is it still possible to construct a valid confidence interval for $\beta_i$ in this case? Can one construct an adaptive honest confidence interval not depending on $k$?

The goal of the present paper is to address these and other related questions on confidence intervals for high-dimensional linear regression with random design. More specifically, we consider construction of confidence intervals for a linear functional $\Tf=\linfun$, where the loading vector $\xi \in \R^{p}$ is given  and $\frac{\max_{i\in \supp(\xi)} |\xi_i|}{\min_{i\in \supp(\xi)} |\xi_i|}\leq \bar{c}$ with $\bar{c}\ge1$ being a constant. Based on the sparsity of $\xi$, we focus on two specific regimes:
the sparse loading regime where $\|\xi\|_0\leq C k$,
with $C>0$ being a constant; the dense loading regime where $\|\xi\|_0$ satisfying \eqref{eq: def of dense support} in Section \ref{sec:formulation}.
It will be seen later that for confidence intervals  $\Tf=\beta_i$ is a prototypical case for the general functional $\Tf=\linfun$ with a sparse loading $\xi$, and $\Tf=\sum_{i=1}^{p} \beta_i$ is a representative case for $\Tf=\linfun$ with a dense loading $\xi$.

To illustrate the main idea, let us first focus on the  two specific functionals $\Tf=\beta_i$ and $\Tf=\sum_{i=1}^{p} \beta_i$.  We establish the convergence rate of the minimax expected length for confidence intervals in the oracle setting where the sparsity parameter $k$ is given.
It is shown that in this case the minimax expected length  is of order $\frac{1}{\sqrt{n}}+k\frac{\log p}{n}$ for confidence intervals for $\beta_i$.  An honest confidence interval, which depends on the sparsity $k$, is constructed and is shown to be minimax rate optimal.
To the best of our knowledge, this is the first construction of confidence intervals in the moderate-sparse region $\msparse$. If the sparsity $k$ falls into the ultra-sparse region $k \lesssim {\sqrt{n}\over \log p}$, the constructed confidence interval is similar to the confidence intervals constructed in \citep{zhang2014confidence,javanmard2014confidence,van2014asymptotically}.
On the other hand, the convergence rate of the minimax expected length of honest confidence intervals for $\sum_{i=1}^{p} \beta_i$ in the oracle setting is shown to be $k\sqrt{\frac{\log p}{n}}$. A rate-optimal confidence interval that also depends on $k$ is constructed. It should be noted that this confidence interval is not based on the de-biased estimator.

One drawback of the constructed confidence intervals mentioned above is that they require prior knowledge of the sparsity $k$. Such knowledge of sparsity is usually unavailable in applications.
A natural question is: Without knowing the sparsity $k$, is it possible to construct a confidence interval as good as when the sparsity $k$ is known? This is a question about adaptive inference, which has been a major goal in nonparametric and high-dimensional statistics. Ideally, an adaptive confidence interval should have its length automatically adjusted to the true sparsity of the unknown regression vector, while maintaining a prespecified coverage probability.  We show that, unlike point estimation,  such a goal is in general not attainable for confidence intervals. In the case of confidence intervals for $\beta_i$, it is impossible to adapt between different sparsity levels,  except when the sparsity $k$ is restricted to the ultra-sparse region $k \lesssim {\sqrt{n}\over \log p}$, over which the confidence intervals have the optimal length of the parametric rate $\frac{1}{\sqrt{n}}$, which does not depend on $k$. In the case of confidence intervals for $\sum_{i=1}^p \beta_i$, it is shown that adaptation to the sparsity is not possible at all, even in the ultra-sparse region $k \lesssim {\sqrt{n}\over \log p}$.

Minimax theory is often criticized as being too conservative as it focuses on the worst case performance. For confidence intervals for high dimensional linear regression, we establish strong non-adaptivity results which demonstrate that the lack of adaptivity is not due to the conservativeness of the minimax framework.  It shows that for any confidence interval with guaranteed coverage probability over the set of $k$ sparse vectors, its expected length at any given point in a large subset of the parameter space must be at least of the same order as the minimax expected length. So the confidence interval must be long at a large subset of points in the parameter space, not just at a small number of  ``unlucky" points. This leads directly to the impossibility of adaptation over different sparsity levels. Fundamentally, the lack of adaptivity is caused by the difficulty in accurately learning the bias of any estimator for high-dimensional linear regression.

We now turn to confidence intervals for general linear functionals.
For a linear functional $\xi^{\intercal} \beta$ in the sparse loading regime,
the rate of the minimax expected length is  $\|\xi\|_2\left(\frac{1}{\sqrt{n}}+k\frac{\log p}{n}\right)$, where $\|\xi\|_2$ is the vector $\ell_2$ norm of $\xi$.
For a linear functional $\xi^{\intercal} \beta$ in the dense loading regime,
the rate of the minimax expected length is  $\|\xi\|_{\infty} k\sqrt{\frac{\log p}{n}}$, where $\|\xi\|_{\infty}$ is the vector $\ell_{\infty}$ norm of $\xi$.
Regarding adaptivity, the phenomena observed in confidence intervals for the two special linear functionals $\Tf=\beta_i$ and $\Tf=\sum_{i=1}^{p} \beta_i$ extend to the general linear functionals. The case of  confidence intervals for $\Tf=\sum_{i=1}^{p}\xi_i \beta_i$ with a sparse loading $\xi$ is similar to that of confidence intervals for $\beta_i$ in the sense that rate-optimal adaptation is impossible except when the sparsity $k$ is restricted to the ultra-sparse region $k\lesssim \frac{\sqrt{n}}{\log p}$.  On the other hand, the case for a dense loading $\xi$ is similar to that of confidence intervals for $\sum_{i=1}^p \beta_i$: adaptation to the sparsity $k$ is not possible at all, even in the ultra-sparse region $k\lesssim \frac{\sqrt{n}}{\log p}$.

In addition to the more typical setting in practice where the covariance matrix $\Sigma$ of the random design and the noise level $\sigma$ of the linear model are unknown, we also consider the case with the prior knowledge of $\Sigma={\rm I}$ and $\sigma=\sigma_0$. It turns out that this case is strikingly different. The minimax rate for the expected length in the sparse loading regime is reduced from $\|\xi\|_2\left(\frac{1}{\sqrt{n}}+k\frac{\log p}{n}\right)$ to $\frac{\|\xi\|_2}{\sqrt{n}}$, and in particular it does not depend on the sparsity $k$. Furthermore, in marked contrast to the case of unknown $\Sigma$ and $\sigma$, adaptation to sparsity  is possible over the full range $k\lesssim \frac{{n}}{\log p}$. On the other hand, for linear functionals $\xi^{\intercal} \beta$  with a dense loading $\xi$, the minimax rates and impossibility for adaptive confidence intervals do not change even with the prior knowledge of $\Sigma={\rm I}$ and $\sigma=\sigma_0$. However, the cost of adaptation is reduced with the prior knowledge.

The rest of the paper is organized as follows: After basic notation is introduced, Section \ref{sec:formulation} presents a precise formulation for the adaptive confidence interval problem. Section \ref{sec:sparse loading regime} establishes the minimaxity and adaptivity results for a general linear functional $\linfun$ with a sparse loading $\xi$. Section \ref{sec:dense loading regime} focuses on confidence intervals for a general linear functional $\linfun$ with a dense loading $\xi$. Section \ref{sec: known design confidence interval} considers the case when there is prior knowledge of covariance matrix of the random design and the noise level of the linear model. Section \ref{sec:discussion} discusses connections to other work and further research directions. The proofs of the main results are given in Section \ref{sec:proof}. More discussion and proofs are presented in the supplement \citep{cai2015supplement}.


\section{Formulation for adaptive confidence interval problem}
\label{sec:formulation}

We present in this section the framework for studying the adaptivity of confidence intervals. We begin with the notation that will be used throughout the paper.

\subsection{Notation}
\label{sec:notation}
For a matrix $X\in \R^{n\times p}$, $X_{i\cdot}$,  $X_{\cdot j}$, and $X_{i,j}$ denote respectively the $i$-th row,  $j$-th column, and  $(i,j)$ entry of the matrix $X$, $X_{i,-j}$ denotes the $i$-th row of $X$ excluding the $j$-th coordinate, and $X_{-j}$ denotes the submatrix of $X$ excluding the $j$-th column. Let $[p]=\{1,2,\cdots,p\}$. For a subset $J\subset[p]$,  $X_{J}$ denotes the submatrix of $X$ consisting of columns  $X_{\cdot j}$ with $j\in J$ and for a vector $x\in \R^{p}$, $x_{J}$ is the subvector of $x$ with indices in $J$ and $x_{-J}$ is the subvector with indices in $J^{c}$.  For a set $S$, $\left|S\right|$ denotes the cardinality of $S$. For a vector $x\in \R^{p}$, $\suppx$ denotes the support of $x$ and the $\ell_q$ norm of $x$ is defined as $\|x\|_{q}=\left(\sum_{i=1}^{q}|x_i|^q\right)^{\frac{1}{q}}$ for $q \geq 0$ with $\|x\|_0=|\suppx|$ and $\|x\|_{\infty}=\max_{1\leq j \leq p}|x_j|$. We use $e_i$ to denote the $i$-th standard basis vector in $\R^p$.
For $a\in \R$, $a_{+}=\max\left\{a,0\right\}$. 
We use $\sum \beta_i$ as a shorthand for $\sum_{i=1}^{p} \beta_i$, $\max\|\Xj\|_2$ as a shorthand for $\max_{1\leq j \leq p}\|\Xj\|_2$ and $\min\|\Xj\|_2$ as a shorthand for $\min_{1\leq j \leq p}\|\Xj\|_2$. For a matrix $A$ and $1\le q\le \infty$, $\|A\|_{q}=\sup_{\|x\|_q = 1} \|Ax\|_q$ is the matrix $\ell_q$ operator norm. In particular, $\|A\|_2$ is the spectral norm. For a symmetric matrix $A$, $\lambda_{\min}\left(A\right)$ and $\lambda_{\max}\left(A\right)$  denote respectively the smallest and largest eigenvalue of $A$.  We use $c$ and $C$ to denote generic positive constants that may vary from place to place. For two positive sequences $a_n$ and $b_n$,  $a_n \lesssim b_n$ means $a_n \leq C b_n$ for all $n$ and $a_n \gtrsim b_n $ if $b_n\lesssim  a_n$ and $a_n \asymp b_n $ if $a_n \lesssim b_n$ and $b_n \lesssim a_n$, and $a_n \ll b_n$ if $\limsup_{n\rightarrow\infty} \frac{a_n}{b_n}=0$ and $a_n \gg b_n$ if $b_n \ll a_n$.

\subsection{Framework for adaptivity of confidence intervals}
\label{sec:framework}
We shall focus in this paper on the high-dimensional linear model with the Gaussian design,
\begin{equation}
y_{n\times 1}=X_{n\times p}\beta_{p\times 1}+\epsilon_{n\times 1}, \quad \epsilon \sim N_n(0,\sigma^2 {\rm I}),
\label{eq: linear model}
\end{equation}
where the rows of $X$ satisfy $X_{i \cdot} \stackrel{\rm i.i.d.}{\sim} N_p(0,\Sigma)$, $i=1, ..., n$, and are independent of $\epsilon$. Both $\Sigma$ and the noise level $\sigma$ are unknown. Let $\Omega=\Sigma^{-1}$ denote the precision matrix. The parameter $\theta=\left(\beta,\Omega, \sigma\right)$ consists of the signal $\beta$, the precision matrix $\Omega$ for the random design, and the noise level $\sigma$. The target of interest is the linear functional of $\beta$, $\T \left(\beta\right)=\xi^{\intercal} \beta,$ where $\xi \in \R^{p}$ is a pre-specified loading vector. The data that we observe is  $Z=\left( Z_1,\cdots, Z_{n}\right)^{\intercal},$ where
$Z_i=\left(y_i,X_i\right)\in \R^{p+1}$ for $i=1,\cdots,n$.

For $0<\alpha<1$ and a given parameter space $\Theta$ and the linear functional $\Tf$, denote by  $\SCIZ$ the set of all $(1-\alpha)$ level confidence intervals for $\T\left(\beta\right)$ over the parameter space $\Theta$, 
\begin{equation}
\SCIZ=\left\{\CIZg=[l(Z), u(Z)]: \inf_{\theta\in \Theta} \PP(l(Z) \le \T(\beta) \le u(Z))\geq 1-\alpha \right\}.
\end{equation}
For any confidence interval $\CIZg \in \SCIZ$, the maximum expected length over a parameter space $\Theta$ is defined as
$$L({\CIZg},\Theta,{\rm T})=\sup_{\theta \in \Theta} \E_{\theta} L\left(\CIZg\right),$$
where for confidence interval $\CIZg=[l(Z), u(Z)]$, $L(\CIZg)= u(Z)-l(Z)$ denotes its length. For two parameter spaces $\Theta_1 \subseteq \Theta$, we define the benchmark $L_{\alpha}^{*}(\Theta_1,\Theta, \Tb)$ as the infimum of the maximum expected length over $\Theta_1$ among all $(1-\alpha)$-level confidence intervals over $\Theta$,
\begin{equation}
L_{\alpha}^{*}(\Theta_1,\Theta,\Tb)=\inf_{\CIZg\in \SCIZ} L(\CIZg, \Theta_1,\Tb).
\label{eq: benchmark}
\end{equation}
We will write $L_{\alpha}^{*}(\Theta,\Tb)$ for $L_{\alpha}^{*}(\Theta,\Theta,\Tb)$, which is the minimax expected length of confidence intervals over $\Theta$. 

We should emphasize that $L_{\alpha}^{*}(\Theta_1,\Theta,\Tb)$ is an important quantity that measures the degree of adaptivity over the nested spaces $\Theta_1\subset \Theta$. A confidence interval $\CIZg$ that is (rate-optimally)  adaptive over $\Theta_1$ and $\Theta$ should have the optimal expected length performance simultaneously over both $\Theta_1$ and $\Theta$ while maintaining a given coverage probability over $\Theta$, i.e.,  $\CIZg\in \SCIZ$ such that
$$L({\CIZg},\Theta_1,\Tb) \asymp L_{\alpha}^{*}(\Theta_1,\Tb) \quad  \text{and} \quad L({\CIZg},\Theta,\Tb) \asymp L_{\alpha}^{*}(\Theta,\Tb).$$ Note that in this case $L({\CIZg},\Theta_1,\Tb)\ge L_{\alpha}^{*}(\Theta_1,\Theta,\Tb)$. So for  two parameter spaces $\Theta_1\subset \Theta$, if $L_{\alpha}^{*}(\Theta_1,\Theta,\Tb) \gg  L_{\alpha}^{*}(\Theta_1,\Tb),
$
then rate-optimal adaptation between $\Theta_1$ and $\Theta$ is impossible to achieve.

We consider the following collection of parameter spaces,
\begin{equation}
\Theta(k)=\left\{\theta=(\beta,\Omega,\sigma):
 \|\beta\|_0\leq k, \frac{1}{\Md}\leq \lambda_{\min} (\Omega)\leq \lambda_{\max}(\Omega) \leq \Md,  0<\sigma\leq \Mv\right\},
 \label{eq: parameter space}
\end{equation}
where $\Md >1$ and $\Mv>0$ are positive constants. Basically, $\Theta(k)$ is the set of all $k$-sparse regression vectors. $\frac{1}{\Md}\leq \lambda_{\min} (\Omega)\leq \lambda_{\max}(\Omega) \leq \Md$ and $0<\sigma\leq \Mv$ are two mild regularity conditions on the design and noise level.

The main goal of this paper is to address the following two questions:
\begin{enumerate}
\item {\it What is the minimax length $L_{\alpha}^{*}(\Theta(k),\Tb)$ in the oracle setting where the sparsity level $k$ is known?} 
\item {\it Is it possible to achieve rate-optimal adaptation over different sparsity levels?} \\
More specifically, for $k_1\ll k$, is it possible to construct a confidence interval $\CIZg$ that is adaptive over $\Theta(k_1)$ and $\Theta(k)$ in the sense that $\CIZg \in \SCIZg$ and
\end{enumerate}
\begin{equation}
\begin{aligned}
L({\CIZg},\Theta(k_1),\Tb)& \asymp L_{\alpha}^{*}(\Theta(k_1),\Tb), \\
 L({\CIZg},\Theta(k),\Tb)& \asymp L_{\alpha}^{*}(\Theta(k),\Tb)?
 \end{aligned}
\label{eq: adaptive benchmark}
\end{equation}
We will answer these questions by analyzing the two benchmark quantities $L_{\alpha}^{*}(\Theta(k),\Tb)$ and $L_{\alpha}^{*}(\Theta(k_1), \Theta(k),\Tb)$. Both lower and upper bounds will be established.
If $\eqref{eq: adaptive benchmark}$ can be achieved, it means that the confidence interval $\CIZg$ can automatically adjust its length to the sparsity level of the true regression vector $\beta$. On the other hand, if $L_{\alpha}^{*}(\Theta(k_1), \Theta(k),\Tb) \gg L_{\alpha}^{*}(\Theta(k_1),\Tb)$, then such a goal is not attainable.

For ease of presentation, we calibrate the sparsity level 
\[
k\asymp p^{\gamma}\quad \mbox{for some $0\leq \gamma<\frac{1}{2}$},
\]
and restrict the loading $\xi$ to the set
$$\xi\in \Xi\left(q,\bar{c}\right)=\left\{\xi \in \R^{p}: \|\xi\|_0= q,  \; \xi \neq {\bf{0}} \; \text{ and } \; \frac{\max_{j\in {\rm supp\left(\xi\right)}}|\xi_j|}{\min_{j\in  {\rm supp\left(\xi\right)}}|\xi_j|}\leq \bar{c}\right\},$$
where $\bar{c}\ge 1$ is a constant. The minimax rate and adaptivity of confidence intervals for the general linear functional $\xi^{\intercal} \beta$ also depends on the sparsity of $\xi$. We are particularly interested in the following two regimes:
\begin{enumerate}
\item The sparse loading regime: $\xi\in \Xi\left(q,\bar{c}\right)$ with
\begin{equation}
q \leq C k.
\label{eq: def of sparse support}
\end{equation}
\item The dense loading regime: $\xi\in \Xi\left(q,\bar{c}\right)$ with
\begin{equation}
q= c p^{\gamma_{q}} \quad \text{with} \quad  2\gamma<\gamma_{q}\leq 1.
\label{eq: def of dense support}
\end{equation}
 \end{enumerate}
The behavior of the problem is significantly different in these two regimes. We will consider separately the sparse loading regime in Section \ref{sec:sparse loading regime} and the dense loading regime in Section \ref{sec:dense loading regime}.

\section{Minimax rate and adaptivity of confidence intervals for sparse loading linear functionals}
\label{sec:sparse loading regime}

In this section, we establish the rates of convergence for the minimax expected length of confidence intervals for $\xi^{\intercal} \beta$ with a sparse loading $\xi$ in the oracle setting where the sparsity parameter $k$ of the regression vector $\beta$ is given. Both minimax upper and lower bounds are given. Confidence intervals for $\xi^{\intercal} \beta$ are constructed and shown to be minimax rate-optimal in the sparse loading regime. Finally, we establish the possibility of adaptivity for the linear functional $\xi^{\intercal} \beta$ with a  sparse loading $\xi$.

\subsection{Minimax length of confidence intervals for $\xi^{\intercal} \beta$ in the sparse loading regime}
\label{sec: minimax sparse}
In this section, we focus on the sparse loading regime defined in \eqref{eq: def of sparse support}. The following theorem establishes the minimax rates for the expected length of confidence intervals for $\xi^{\intercal} \beta$ in the sparse loading regime.

\begin{Theorem}
\label{thm: minimax for sparse linear functional}
Suppose that $0<\alpha<\frac{1}{2}$ and $k \leq c \min\{p^{\gamma},\frac{n}{\log p}\}$ for some constants $c>0$ and $ 0\leq \gamma < \frac{1}{2}$. If $\xi$ belongs to the sparse loading regime \eqref{eq: def of sparse support}, the minimax expected length for $\left(1-\alpha\right)$ level confidence intervals of $\xi^{\intercal} \beta$ over $\Theta\left(k\right)$ satisfies
\begin{equation}
L_{\alpha}^{*}\left(\Theta\left({k}\right),\xi^{\intercal} \beta\right)\asymp \|\xi\|_2 \left( \frac{1}{\sqrt{n}}+k\frac{\log p}{n}\right).
\label{eq: minimax for multi sparse}
\end{equation}
\end{Theorem}
Theorem \ref{thm: minimax for sparse linear functional} is established in two separate steps.
\begin{enumerate}
\item Minimax upper bound: we construct a confidence interval $\CIZs$  such that $\CIZs \in \SCIZl$ and for some constant $C>0$
\begin{equation}
L\left(\CIZs, \Theta\left(k\right),\xi^{\intercal} \beta\right)\leq C \|\xi\|_2 \left(\frac{1}{\sqrt{n}}+k\frac{\log p}{n}\right).
\end{equation}
\item Minimax lower bound: we show that for some constant $c>0$
\begin{equation}
L_{\alpha}^{*}\left(\Theta\left(k\right),\xi^{\intercal} \beta \right) \geq c \|\xi\|_2 \left(\frac{1}{\sqrt{n}}+k\frac{\log p}{n}\right).
\label{eq: lower bound for beta_1}
\end{equation}
\end{enumerate}
The minimax lower bound is implied by the adaptivity result given in Theorem \ref{thm: non-adaptivity multi sparse}. 
We now detail the construction of a confidence interval $\CIZs$ achieving the minimax rate \eqref{eq: minimax for multi sparse} in the sparse loading regime. The  interval $\CIZs$ is centered at a de-biased scaled Lasso estimator, which generalizes the ideas used in \citep{zhang2014confidence,javanmard2014confidence, van2014asymptotically}. The construction of the (random) length is different from the aforementioned papers as the asymptotic normality result is not valid once $k\gtrsim {\sqrt{n}\over \log p}$.

 Let $\{\widehat{\beta},\hat{\sigma}\}$ be the scaled Lasso estimator with $\lambda_0= \sqrt{\frac{2.05 \log p}{n}}$,
\begin{equation}
\{\widehat{\beta},\hat{\sigma}\}=\argmin_{\beta \in \R^{p},\sigma
\in \R^{+}}\frac{\|y-X\beta\|_2^2}{2n\sigma}+\frac{\sigma}{2}+\lambda_0 \sum_{j=1}^{p} \frac{\|\Xj\|_2}{\sqrt{n}} |\beta_j|.
\label{eq: Scaled Lasso}
\end{equation}
Define
\begin{equation}
\begin{aligned}
\widehat{\uu}=&\argmin_{\uu\in \R^p} \left\{\uu^{\intercal} \widehat{\Sigma} \uu: \|\widehat{\Sigma}\uu-\xi\|_{\infty}\leq \lambda_n\right\},
\end{aligned}
\label{eq: DZ for linear functional}
\end{equation}
where $\widehat{\Sigma}=\frac{1}{n} X^{\intercal} X$ and $\lambda_n=12 \|\xi\|_2 \Md^2 \sqrt{\frac{\log p}{n}}$.
The confidence interval $\CIZs$ is centered at the following de-biased estimator
\begin{equation}
\widetilde{\mu}=\xi^{\intercal} \widehat{\beta}+\widehat{\uu}^{\intercal} \frac{1}{n}X^{\intercal}\left(y-X\widehat{\beta}\right),
\label{eq: def of debiased estimator of linear functional}
\end{equation}
where $\widehat{\beta}$ is the scaled Lasso estimator given in \eqref{eq: Scaled Lasso} and $\widehat{\uu}$ is defined in \eqref{eq: DZ for linear functional}.
Before specifying the length of the confidence interval, we review the following definition of restricted eigenvalue introduced in \citep{bickel2009simultaneous}, 
\begin{equation}
\kappa(X,k,\alpha_0)=\min_{\substack{J_{0}\subset \{1,\cdots,p\},\\|J_0|\leq k}} \min_{\substack{\delta\neq 0,\\ \|\delta_{J_0^c}\|_1\leq \alpha_0 \|\delta_{J_0}\|_1}} \frac{\|X\delta\|_2}{\sqrt{n} \|\delta_{J_{0}}\|_2}.
\end{equation}
Define
\begin{equation}
\rho_1\left(k\right)=\|\xi\|_2 \hat{\sigma} \min \left\{1.01\sqrt{\frac{\widehat{\uu}^{\intercal} \widehat{\Sigma}\widehat{\uu}}{{n\|\xi\|_2^2}}} \za+C_1\left(X,k\right) k \frac{\log p}{n} , \;  \log p (\frac{1}{\sqrt{n}}+\frac{k \log p}{n} ) \right\},
\label{eq: CI length multi sparse}
\end{equation}
where $\za$ is the $\alpha/2$ upper quantile of the standard normal distribution and
\begin{equation}
\label{eq: def of constant for beta_1}
C_1\left(X,k\right)=7000\Md^2\frac{\sqrt{n}}{\min\|\Xj\|_2} \max\left\{1.25,\frac{912\max \|\Xj\|_2^2}{n \kappa^2\left(X,k,405\left(\frac{\max \|\Xj\|_2}{\min \|\Xj\|_2}\right)\right) }\right\}.
\end{equation}
Define the event \begin{equation}
A=\left\{\hat{\sigma}\leq \log p\right\}.
\label{eq: def of event A}
\end{equation}
The confidence interval $\CIZs$ for $\xi^{\intercal} \beta$ is defined as
\begin{equation}
\CIZs= \left\{
  \begin{array}{cl}
    \left[\widetilde{\mu}-\rho_1\left(k\right), \quad \widetilde{\mu}+\rho_1\left(k\right) \right]    &\quad \text{on}\; A\\
   \left\{0\right\}  &\quad \text{on}\; A^{c}
  \end{array}
\right.
\label{eq: constructed CI for multi sparse}
\end{equation}
It will be shown in Section \ref{sec:proof} that the confidence interval $\CIZs$ has the desired coverage property and achieves the minimax length in \eqref{eq: minimax for multi sparse}.

\begin{Remark}\rm
In the special case of $\xi=e_1$, the confidence interval defined in \eqref{eq: constructed CI for multi sparse} is similar to the ones based on the de-biased estimators introduced in \citep{zhang2014confidence, javanmard2014confidence, van2014asymptotically}.
The second term $\widehat{\uu}^{\intercal}\frac{1}{n}X^{\intercal}\left(y-X\widehat{\beta}\right)$ in \eqref{eq: def of debiased estimator of linear functional}  is incorporated to reduce the bias of the scaled Lasso estimator $\widehat{\beta}$.
The constrained estimator $\widehat{\uu}$ defined in \eqref{eq: DZ for linear functional} is a score vector $\uu$ such that the variance term $\uu^{\intercal} \widehat{\Sigma} \uu$ is minimized and one component of the bias term $\|\widehat{\Sigma}\uu-\xi\|_{\infty}$ is constrained by the tuning parameter $\lambda_n$. 
The tuning parameter $\lambda_n$ is chosen as $12 \|\xi\|_2\Md^2 \sqrt{\frac{\log p}{n}}$ such that $\uu=\Omega \xi$ lies in the constraint set  $\|\widehat{\Sigma}\uu-\xi\|_{\infty}\leq \lambda_n$ in \eqref{eq: DZ for linear functional} with overwhelming probability. 
For $C_1(X,k)$ defined in \eqref{eq: def of constant for beta_1}, it will be shown that it is upper bounded by a constant with overwhelming probability.
\end{Remark}
\subsection{Adaptivity of confidence intervals for $\xi^{\intercal} \beta$ in the sparse loading regime}
\label{sec: sparse adaptivity}
We have constructed a minimax rate-optimal confidence interval for $\xi^{\intercal} \beta$ in the oracle setting where the sparsity $k$ is assumed to be known. A major drawback of the construction is that it requires prior knowledge of $k$, which is typically unavailable in practice. An interesting question is whether it is possible to construct adaptive confidence intervals that have the guaranteed coverage and automatically adjust its length to $k$.

We now consider the adaptivity of the confidence intervals for $\xi^{\intercal} \beta$.
In light of the minimax expected length given in Theorem \ref{thm: minimax for sparse linear functional},
the following theorem provides an answer to the adaptivity question \eqref{eq: adaptive benchmark} for the confidence intervals for $\xi^{\intercal}\beta$ in the sparse loading regime.

\begin{figure}[htb]
 \centering
    \includegraphics[width=4.5in,height=1.75in]{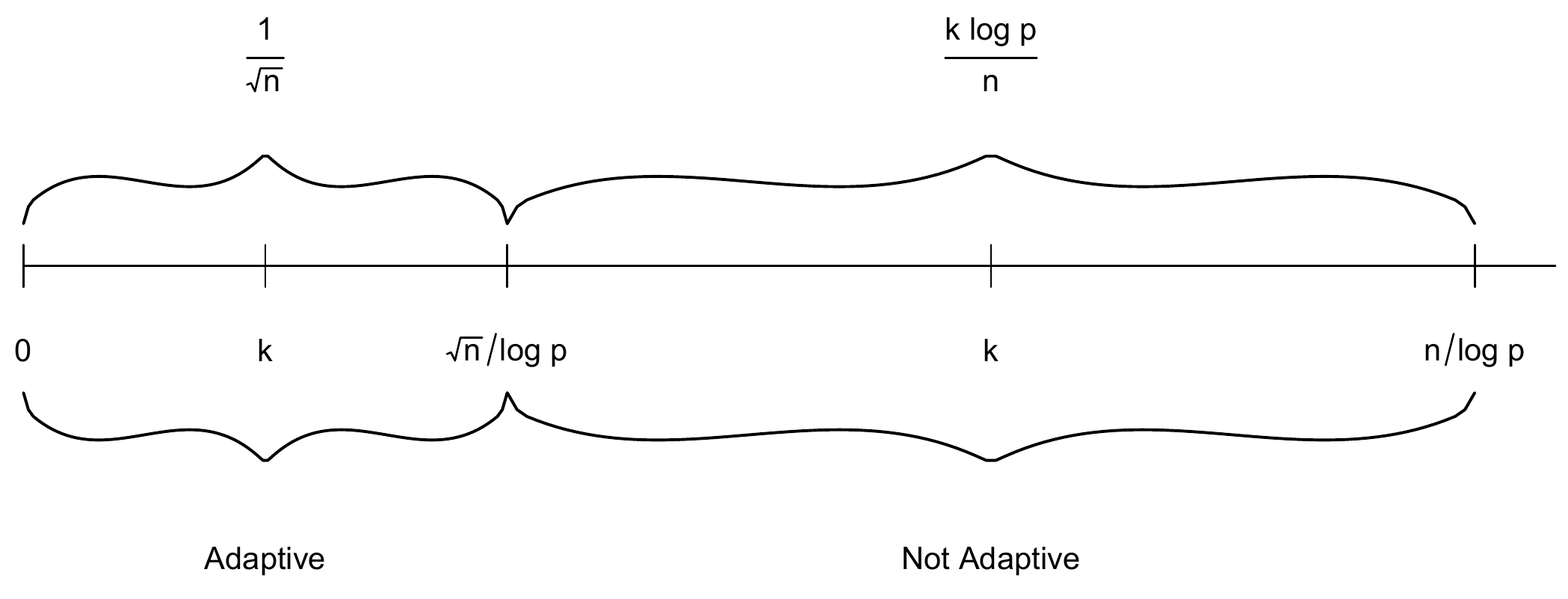}
    \caption{Illustration of adaptivity of confidence intervals for $\xi^{\intercal}\beta$ with a sparse loading $\xi$. For adaptation between  $\Theta(k_1)$ and $\Theta(k)$ with $k_1\ll k$, rate-optimal adaptation is possible if $\usparse$ and impossible otherwise.}
    \label{fig: adaptivityplot}
\end{figure}

\begin{Theorem}
\label{thm: non-adaptivity multi sparse}
Suppose that $0 < \alpha < \frac{1}{2}$ and $k_1\le k \leq c\min\left\{p^{\gamma},\frac{n}{\log p}\right\}$ for some constants $c>0$ and $0\leq \gamma <\frac{1}{2}$. Then
\begin{equation}
L_{\alpha}^{*}(\Theta(k_1),\Theta(k),\xi^{\intercal} \beta)\geq c_1 \|\xi\|_2\left( \frac{1}{\sqrt{n}}+k\frac{\log p}{n}\right),
\label{eq: non-adaptivity multi sparse}
\end{equation}
for some constant $c_1>0$.
\end{Theorem}
Note that Theorem \ref{thm: non-adaptivity multi sparse} implies the minimax lower bound in Theorem \ref{thm: minimax for sparse linear functional} by taking $k_1=k$.
Theorem \ref{thm: non-adaptivity multi sparse} rules out the possibility of rate-optimal adaptive confidence intervals beyond the ultra-sparse region. Consider the setting where $k_1\ll k$ and $\msparse$. In this case,
\[
L_{\alpha}^{*}(\Theta(k_1),\Theta(k),\xi^{\intercal}\beta) \asymp L_{\alpha}^{*}(\Theta(k),\xi^{\intercal}\beta)  \asymp  \|\xi\|_2 k \frac{\log p}{n}  \gg  L_{\alpha}^{*}(\Theta(k_1),\xi^{\intercal}\beta).
\]
So it is impossible to construct a confidence interval that is adaptive simultaneously over $\Theta(k_1)$ and $\Theta(k)$ when $\msparse$ and  $k_1\ll k$. The only possible region for adaptation is in the ultra-sparse region $\usparse$, over which the optimal expected length of confidence intervals is of order $1\over \sqrt{n}$ and in particular does not depend on the specific sparsity level. These facts are illustrated in Figure \ref{fig: adaptivityplot}.

So far the analysis is carried out within the minimax framework where the focus is on the performance in the worst case over a large parameter space. The minimax theory is often criticized as being too conservative. In the following, we establish a stronger version of the non-adaptivity result which demonstrates that the lack of adaptivity for confidence intervals is not due to the conservativeness of the minimax framework. The result shows that for any confidence interval $\CIZl$, under the coverage constraint that $\CIZl \in \SCIZl$, its expected length at any given $\theta^{*}=\left(\beta^*,\rm I, \sigma \right) \in  \Theta\left(k_1\right)$ must be of order
$\|\xi\|_2\left(\frac{1}{\sqrt{n}}  +  {k} \frac{\log p}{n}\right).$ So the confidence interval must be long at a large subset of points in the parameter space, not just at a small number of ``unlucky" points.

\begin{Theorem}
\label{thm: strong non-adaptivity multi sparse}
Suppose that $0<\alpha<\frac{1}{2}$ and $k \leq c\min\{p^{\gamma},\frac{n}{\log p}\}$ for some constants $c>0$ and $ 0\leq \gamma< \frac{1}{2}$. Let $k_1\leq \left(1-\zeta_0\right)k-1$ and $q\leq \frac{\zeta_0}{4} k$ for some constant $0<\zeta_0<1$. Then  for any  $\theta^{*}=\left(\beta^{*},\rm I, \sigma \right) \in  \Theta\left(k_1\right)$ and $\xi \in \Xi\left(q,\bar{c}\right)$,
\begin{equation}
\inf_{\CIZl \in \SCIZl}\E_{\theta^*} L\left(\CIZl \right) \geq c_1 \|\xi\|_{2}  \left( {k} \frac{\log p}{n}+\frac{1}{\sqrt{n}} \right)\sigma,
\label{eq: strong non-adaptivity multi sparse}
\end{equation}
for some constant $c_1>0$.
\end{Theorem}
Note that no supremum is taken over the parameter $\theta^*$ in \eqref{eq: strong non-adaptivity multi sparse}. Theorem \ref{thm: strong non-adaptivity multi sparse} illustrates that if a confidence interval $\CIZl$ is ``superefficient" at any point $\theta^{*}=\left(\beta^*,\rm I, \sigma \right) \in \Theta(k_1)$ in the sense that
$$\E_{\theta^*} L\left(\CIZl\right)\ll  \|\xi\|_2\left( \frac{1}{\sqrt{n}}+k\frac{\log p}{n}\right)\sigma,$$
then the confidence interval $\CIZl$ can not have the guaranteed coverage over the parameter space $\Theta(k)$.


\subsection{Minimax rate and adaptivity of confidence intervals for $\beta_1$}
\label{sec: minimax and adaptivity beta1}
We now turn to the special case $\Tf=\beta_i$, which has been the focus of several previous papers  \citep{zhang2014confidence,javanmard2013hypothesis,javanmard2014confidence,van2014asymptotically}. Without loss of generality, we consider $\beta_1$, the first coordinate of $\beta$,  in the following discussion and the results for any other coordinate $\beta_i$ are the same.  The linear functional $\beta_1$ is the special case of linear functional of sparse loading regime with $\xi=e_1$.

Theorem \ref{thm: minimax for sparse linear functional} implies that the minimax expected length for $\left(1-\alpha\right)$ level confidence intervals of $\beta_1$ over $\Theta\left(k\right)$ satisfies\begin{equation}
L_{\alpha}^{*}\left(\Theta\left(k\right),\beta_1\right) \asymp \frac{1}{\sqrt{n}}+k\frac{\log p}{n}.
\label{eq: minimax for first coef in regression}
\end{equation}
In the ultra-sparse region with  $\usparse$, the minimax expected length is of order $\frac{1}{\sqrt{n}}$. However, when $k$ falls in the moderate-sparse region $\msparse$, the minimax expected length is of order $k \frac{\log p}{n}$ and in this case  $k \frac{\log p}{n} \gg \frac{1}{\sqrt{n}}$. Hence the confidence intervals constructed in \citep{zhang2014confidence,javanmard2013hypothesis,javanmard2014confidence,van2014asymptotically}, which are of parametric length $\frac{1}{\sqrt{n}}$, asymptotically have coverage probability going to 0. The condition $\usparse$ is necessary for the parametric rate $\frac{1}{\sqrt{n}}$.
\citep{van2014asymptotically} established asymptotic normality and asymptotic efficiency for a de-biased estimator under the sparsity assumption $k \ll \frac{\sqrt{n}}{\log p}$. Similar results have also been given in \citep{ren2013asymptotic} for a related problem of estimating  a single entry of a $p$-dimensional precision matrix based on $n$ i.i.d. samples under the same sparsity condition $k\ll \frac{\sqrt{n}}{\log p}$. It was also shown that $k \ll \frac{\sqrt{n}}{\log p}$ is necessary for the asymptotic normality and asymptotic efficiency results.

The following corollary, as a special case of Theorem \ref{thm: strong non-adaptivity multi sparse}, illustrates the strong non-adaptivity  for confidence intervals of $\beta_1$ when $k\gg {\sqrt{n}\over \log p}$.
 \begin{Corollary}
 \label{cor: strong non-adaptivity for beta_1}
Suppose that $0<\alpha < \frac{1}{2}$ and  $k \leq c \min\{p^{\gamma},\frac{n}{\log p}\}$ for some constants $c>0$ and $ 0\leq \gamma< \frac{1}{2}$. Let $k_1\leq \left(1-\zeta_0\right)k-1$ for some constant $0<\zeta_0<1$.
Then for any $\theta^{*}=\left(\beta^*,{\rm I}, \sigma \right) \in  \Theta\left(k_1\right)$,
\begin{equation}
\inf_{\CIZa \in \SCIZa }\E_{\theta^*} L\left(\CIZa \right) \geq c_1 \left( \frac{1}{\sqrt{n}}  + {k} \frac{\log p}{n}\right)\sigma,
\label{eq: strong non-adaptivity for beta_1}
\end{equation}
for some constant $c_1>0.$
\end{Corollary}

\section{Minimax rate and adaptivity of confidence intervals for dense loading linear functionals}
\label{sec:dense loading regime}


We now turn to the setting where the loading $\xi$ is  dense in the sense of \eqref{eq: def of dense support}. We will also briefly discuss the special case $\sum_{i=1}^{p} \beta_i$ and the computationally feasible confidence intervals.
\subsection{Minimax length of confidence intervals for $\xi^{\intercal} \beta$ in the dense loading regime}
\label{sec: minimax dense}
The following theorem establishes the minimax length of confidence intervals of $\xi^{\intercal} \beta$ in the dense loading regime \eqref{eq: def of dense support}.
\begin{Theorem}
\label{thm: minimax for dense linear functional}
Suppose that $0<\alpha<\frac{1}{2}$ and $k \leq c \min\{p^{\gamma},\frac{n}{\log p}\}$ for some constants $c>0$ and $ 0\leq \gamma < \frac{1}{2}$. If $\xi$ belongs to the dense loading regime \eqref{eq: def of dense support}, the minimax expected length for $\left(1-\alpha\right)$ level confidence intervals of $\xi^{\intercal} \beta$ over $\Theta\left(k\right)$ satisfies
\begin{equation}
L_{\alpha}^{*}\left(\Theta\left(k\right),\xi^{\intercal} \beta\right)\asymp \|\xi\|_\infty  k\sqrt{\frac{\log p}{n}}.
\label{eq: minimax for multi dense}
\end{equation}
\end{Theorem}
Note that the minimax rate in  \eqref{eq: minimax for multi dense} is significantly different from the minimax rate $\|\xi\|_2( \frac{1}{\sqrt{n}}+k\frac{\log p}{n})$ for the sparse loading case given in Theorem \ref{thm: minimax for sparse linear functional}.
In the following, we construct a confidence interval $\CIZd$ achieving the minimax rate \eqref{eq: minimax for multi dense} in the dense loading regime. Define
\begin{equation}
C_2(X,k)=822\frac{\sqrt{n}}{\min\|\Xj\|_2} \max\left\{1.25,\frac{912\max \|\Xj\|_2^2 }{n \kappa^2\left(X,k,405\left(\frac{\max \|\Xj\|_2}{\min \|\Xj\|_2}\right)\right) }\right\}.
\label{eq: def of constant for sum beta}
\end{equation}
It will be shown that $C_2(X,k)$ is upper bounded by a constant with overwhelming probability.
The confidence interval $\CIZd$ is defined to be,
\begin{equation}
\CIZd
=\left\{\begin{array}{cl}
\left[\xi^{\intercal} \widehat{\beta}-\|\xi\|_{\infty}  \rho_2\left(k\right), \xi^{\intercal} \widehat{\beta}+\|\xi\|_{\infty} \rho_2\left(k\right) \right]& \quad \text{on} \; A\\
\left\{0\right\} & \quad \text{on} \; A^c
\end{array}
\right.
\label{eq: constructed CI for multi dense}
\end{equation}
where $A$ is defined in \eqref{eq: def of event A}  and $\widehat{\beta}$ is the scaled Lasso estimator defined in \eqref{eq: Scaled Lasso} and  
\begin{equation}
\rho_2\left(k\right)=\min\left\{C_{2}\left(X,k\right)k\sqrt{\frac{\log p}{n}}\hat{\sigma}, \log p \left(k \sqrt{\frac{\log p}{n}} \hat{\sigma}\right)\right\}.
\label{eq: CI length sum beta}
\end{equation}
The confidence interval constructed in \eqref{eq: constructed CI for multi dense} will be shown to have the desired coverage property and achieve the minimax length in \eqref{eq: minimax for multi dense}.
A major difference between the construction of $\CIZd$ and that of $\CIZs$ is that $\CIZd$ is not centered at a de-biased estimator.
If a de-biased estimator is used for the construction of confidence intervals for $\xi^{\intercal} \beta$ with a dense loading, its variance would be too large, much larger than the optimal length $\|\xi\|_{\infty} k\sqrt{\frac{\log p}{n}}$.
\subsection{Adaptivity of confidence intervals for $\xi^{\intercal} \beta$ in the dense loading regime}
\label{sec: adaptivity dense}
In this section, we investigate the possibility of adaptive confidence intervals for $\xi^{\intercal} \beta$ in the dense loading regime.
The following theorem leads directly to an answer to  the adaptivity question \eqref{eq: adaptive benchmark} for confidence intervals for $\xi^{\intercal} \beta$ in the dense loading regime.

\begin{Theorem}
\label{thm: adaptivity multi dense}
Suppose that $0 < \alpha <\frac{1}{2}$ and $k_1\le k \leq c\min\left\{p^{\gamma},\frac{n}{\log p}\right\}$ for some constants $c>0$ and $0\leq \gamma <\frac{1}{2}$. Then, for some constant $c_1>0$,
\begin{equation}
\label{eq: adaptivity multi dense}
L_{\alpha}^{*}\left(\Theta\left({k_1}\right),\Theta\left({k}\right),\xi^{\intercal} \beta \right)\geq c_1 \|\xi\|_{\infty}{k} \sqrt{\frac{\log p}{n}}.
\end{equation}
\end{Theorem}
Theorem \ref{thm: adaptivity multi dense} implies the minimax lower bound in Theorem \ref{thm: minimax for dense linear functional} by taking $k_1=k$.
If $k_1\ll k,$ \eqref{eq: adaptivity multi dense} implies 
\begin{equation}
L_{\alpha}^{*}\left(\Theta\left({k_1}\right),\Theta\left({k}\right),\xi^{\intercal}\beta\right)\ge c \|\xi\|_{\infty}k\sqrt{\frac{\log p}{n}} \gg  L_{\alpha}^{*}\left(\Theta\left({k_1}\right),\xi^{\intercal}\beta\right),
\label{eq: impossibility of adaptivity}
\end{equation}
which shows that rate-optimal adaptation over two different sparsity levels $k_1$ and $k$ is not possible at all for any $k_1\ll k$.
In contrast, in the case of the sparse loading regime,  Theorem \ref{thm: non-adaptivity multi sparse} shows that it is possible to construct an adaptive confidence interval in the ultra-sparse region $\usparse$, although adaptation is not possible in the moderate-sparse region $\msparse$.

Similarly to Theorem \ref{thm: strong non-adaptivity multi sparse}, the following theorem establishes the strong non-adaptivity results for $\xi^{\intercal} \beta$ in the dense loading regime.
 \begin{Theorem}
 \label{thm: strong non-adaptivity multi dense}
Suppose that $0<\alpha<\frac{1}{2}$ and $k \leq c\min\{p^{\gamma},\frac{n}{\log p}\}$ for some constants $c>0$ and $ 0\leq \gamma< \frac{1}{2}$. Let $q$ satisfies \eqref{eq: def of dense support} and $k_1\leq \left(1-\zeta_0\right)k-1$ for some positive constant $0<\zeta_0<1$. Then for any $\theta^{*}=\left(\beta^{*},\rm I, \sigma \right) \in  \Theta\left(k_1\right)$ and $\xi \in \Xi\left(q,\bar{c}\right)$,
there is some constant $c_1>0$ such that
\begin{equation}
\inf_{\CIZl \in \SCIZl }\E_{\theta^*} L\left(\CIZl \right) \geq c_1 \|\xi\|_{\infty} {k} \sqrt{\frac{\log p}{n}}\sigma.
\label{eq: strong non-adaptivity multi dense}
\end{equation}
\end{Theorem}

\subsection{Minimax length and adaptivity of confidence intervals for $\sum_{i=1}^{p} \beta_i$}
\label{sec: minimax and adaptivity sumbeta}
We now turn to to the special case of  $T(\beta)=\sum_{i=1}^{p} \beta_i$, the sum of all coefficients.
Theorem \ref{thm: minimax for dense linear functional} implies that the minimax expected length for $(1-\alpha)$ level confidence intervals of $\sum_{i=1}^{p} \beta_i$ over $\Theta\left(k\right)$ satisfies
\begin{equation}
L_{\alpha}^{*}\left(\Theta\left(k\right),\sum \beta_i\right) \asymp {k}\sqrt{\frac{\log p}{n}}. 
\label{eq: minimax for sum coef in regression}
\end{equation}
The following impossibility of adaptivity result for confidence intervals for $\sum_{i=1}^{p}\beta_i$ is a special case of Theorem \ref{thm: strong non-adaptivity multi dense}.
 \begin{Corollary}
 \label{cor: strong non-adaptivity for sum coef}
Suppose that $0<\alpha < \frac{1}{2}$ and  $k \leq c \min\{p^{\gamma},\frac{n}{\log p}\}$ for some constants $c>0$ and $ 0\leq \gamma< \frac{1}{2}$. Let $k_1\leq \left(1-\zeta_0\right)k-1$ for some constant $0<\zeta_0<1$.
Then for any $\theta^{*}=\left(\beta^*,{\rm I}, \sigma \right) \in  \Theta\left(k_1\right)$,
\begin{equation}
\inf_{\CIZb \in \SCIZb}\E_{\theta^*} L\left(\CIZb \right) \geq c_1 {k} \sqrt{\frac{\log p}{n}} \sigma,
\label{eq: strong adaptivity for sum coef in regression1}
\end{equation}
for some constant $c_1>0$.
\end{Corollary}
\begin{Remark}\rm
In the Gaussian sequence model, the problem of estimating the sum of sparse means has been considered in \citep{cai2004minimax,cai2005adaptive} and more recently in \citep{collier2015minimax}. In particular, minimax rate is given in \citep{cai2004minimax} and \citep{collier2015minimax}. The problem of constructing minimax confidence intervals for the sum of sparse normal means was studied in \citep{cai2004adaptation}.
\end{Remark}

\subsection{Computationally feasible confidence intervals}
\label{sec: computational CI}

A major drawback of the minimax rate-optimal confidence intervals $\CIZs$ given in \eqref{eq: constructed CI for multi sparse} and $\CIZd$ given in \eqref{eq: constructed CI for multi dense} is that they are not computationally feasible as both depend on  restricted eigenvalue $\kappa(X, k, \alpha_0)$, which is difficult to evaluate. In this section, we assume the prior knowledge of the sparsity $k$ and discuss how to construct a computationally feasible confidence interval. 

The main idea is to replace the term involved with restricted eigenvalue by a computationally feasible lower bound function $\omega\left(\Omega,X,k\right)$ defined by
\begin{equation}
\label{eq: def of omega function}
\omega(\Omega,X,k)=\left(\frac{1}{4\sqrt{\lambda_{\max}\left(\Omega\right)}}-\frac{9\left(1+405 \frac{\max \|\Xj\|_2}{\min \|\Xj\|_2}\right)}{\sqrt{\lambda_{\min}\left(\Omega\right)}} \sqrt{k \frac{\log p}{n}}\right)_{+}^2.
\end{equation}
The lower bound relation is established by Lemma 13 in the supplement \citep{cai2015supplement}, which is based on the concentration inequality for Gaussian design in \citep{raskutti2010restricted}. Except for $\lambda_{\min}\left(\Omega\right)$ and $\lambda_{\max}\left(\Omega\right)$, all terms in \eqref{eq: def of omega function} are based on the data $\left(X,y\right)$ and the prior knowledge of $k$. To construct a data-dependent computationally feasible confidence interval, we make the following assumption, 
\begin{equation}
\sup_{\Omega\in \GG_{\Omega}}\mathbb{P}_{X}\left(\max\left\{\left|\widetilde{\lambda_{\min}\left(\Omega\right)}-\lambda_{\min}\left(\Omega\right)\right|,\left|\widetilde{\lambda_{\max}\left(\Omega\right)}-\lambda_{\max}\left(\Omega\right)\right|\right\}\geq C a_{n,p}\right)=o(1),
\label{eq: assumption for computation}
\end{equation}
where $\limsup a_{n,p}=0$ and $\GG_{\Omega}$ is a pre-specified parameter space for $\Omega$ and $\mathbb{P}_{X}$ denotes the probability distribution with respect to $X$.
\begin{Remark} \rm
We assume $\GG_{\Omega}$ is a subspace of the precision matrix defined in \eqref{eq: parameter space},  $\left\{\Omega: \frac{1}{\Md} \leq \lambda_{\min}\left(\Omega\right) \leq \lambda_{\max}\left(\Omega\right)\leq \Md\right\}$. By assuming $\GG_{\Omega}$ is the set of precision matrix of special structure, we can find estimators satisfying \eqref{eq: assumption for computation}.  If $\GG_{\Omega}$ is assumed to be the set of sparse precision matrix, we can estimate the precision matrix $\Omega$ by CLIME estimator $\widetilde{\Omega}$ proposed in \citep{cai2011constrained}. Under proper sparsity assumption on $\Omega$, the plugin estimator $\left(\lambda_{\min}\left(\widetilde{\Omega}\right),\lambda_{\max}\left(\widetilde{\Omega}\right)\right)$ satisfies \eqref{eq: assumption for computation}.  Other special structures can also be assumed, for example, the covariance matrix is sparse. We can use the plugin estimator of the estimator proposed in \citep{cai2012optimal}.
\end{Remark}
With $\widetilde{\lambda_{\min}\left(\Omega\right)}$ and $\widetilde{\lambda_{\max}\left(\Omega\right)}$, we define $\widetilde{\omega}\left(\Omega,X,k\right)$ as
$$
\widetilde{\omega}(\Omega,X,k)=\left(\frac{1}{4\sqrt{\widetilde{\lambda_{\max}\left(\Omega\right)}}}-\frac{9\left(1+405 \frac{\max \|\Xj\|_2}{\min \|\Xj\|_2}\right)}{\sqrt{\widetilde{\lambda_{\min}\left(\Omega\right)}}} \sqrt{k \frac{\log p}{n}}\right)_{+}^2.
$$
and  construct computationally feasible confidence intervals  by replacing $\kappa^2\left(X,k,405\left(\frac{\max \|\Xj\|_2}{\min \|\Xj\|_2}\right)\right) $ in \eqref{eq: constructed CI for multi sparse} and \eqref{eq: constructed CI for multi dense} with $\widetilde{\omega}(\Omega,X,k)$. 

\section{Confidence intervals for linear functionals with prior knowledge $\Omega= {\rm I}$ and $\sigma=\sigma_0$} 
\label{sec: known design confidence interval}

We have so far focused on the setting where both  the precision matrix $\Omega$ and the noise level $\sigma$ are unknown, which is the case in most statistical applications. It is still of theoretical interest to study the problem when  $\Omega$ and $\sigma$ are known. It is interesting to contrast the results with the ones when $\Omega$ and $\sigma$ are unknown.
In this case, we consider the setting where it is known a priori that $\Omega= {\rm I}$ and $\sigma=\sigma_0$ and specify the parameter space as
\begin{equation}
\Theta(k,{\rm I},\sigma_0)=\{\theta=(\beta,{\rm I},\sigma_0): \|\beta\|_0\leq k \}.
\label{eq: both known parameter space}
\end{equation}
We will discuss separately the minimax rates and adaptivity of confidence intervals for the linear functionals in the sparse loading regime and dense loading regime over the parameter space  $\Theta(k,{\rm I},\sigma_0)$.
\subsection{Confidence intervals for linear functionals in the sparse loading regime}
\label{sec: sparse functional known}
The following theorem establishes the minimax rate of confidence intervals for linear functionals in the sparse loading regime when there is prior knowledge that $\Omega={\rm I}$ and $\sigma=\sigma_0$.
\begin{Theorem}
\label{thm: minimax for both known}
Suppose that $0<\alpha<\frac{1}{2}$ and $k \leq c \min\{p^{\gamma},\frac{n}{\log p}\}$ for some constants $c>0$ and $ 0\leq \gamma < \frac{1}{2}$. If $\xi$ belongs to the sparse loading regime \eqref{eq: def of sparse support}, the minimax expected length for $\left(1-\alpha\right)$ level confidence intervals of $\xi^{\intercal} \beta$ over $\Theta(k,{\rm I},\sigma_0)$ satisfies
\begin{equation}
L_{\alpha}^{*}\left(\Theta(k,{\rm I},\sigma_0),\xi^{\intercal} \beta\right)\asymp  \frac{\|\xi\|_2}{\sqrt{n}}.
\label{eq: minimax for known Sigma}
\end{equation}
\end{Theorem}
Compared with the minimax rate $\frac{\|\xi\|_2}{\sqrt{n}}+\|\xi\|_2k\frac{\log p}{n}$ for the unknown  $\Omega$ and $\sigma$ case given in Theorem \ref{thm: minimax for sparse linear functional}, the minimax rate in  \eqref{eq: minimax for known Sigma} is significantly different.
With the prior knowledge of $\Omega={\rm I}$ and $\sigma=\sigma_0$, the above theorem shows that the minimax expected length of confidence intervals for $\xi^{\intercal} \beta$ is always of parametric rate and in particular does not depend on the sparsity parameter $k$. In this case, adaptive confidence intervals for $\xi^{\intercal} \beta$ is possible over the full range $k\leq c\frac{n}{\log p}$. A similar result for confidence intervals covering all $\beta_i$ has been given in a recent paper \citep{javanmard2015biasing}. The focus of \citep{javanmard2015biasing} is on individual coordinates, not general linear functionals.

The minimax lower bound of Theorem \ref{thm: minimax for both known} follows from the parametric lower bound of Theorem \ref{thm: minimax for sparse linear functional}.  As both  $\Omega$ and $\sigma$ are known, the upper bound analysis is easier than the unknown $\Omega$ and $\sigma$ case and is similar to the one given in  \citep{javanmard2015biasing}.  For completeness, we detail the construction of a confidence interval achieving the minimax length in \eqref{eq: minimax for known Sigma} using the de-biasing method.
We first randomly split the samples $\left(X,y\right)$ into two subsamples $\left(\xa,\ya\right)$ and $\left(\xb,\yb\right)$ with sample sizes $n_1$ and $n_2$, respectively. Without loss of generality, we assume that $n$ is even and $n_1=n_2=\frac{n}{2}$.  Let $\widehat{\beta}$ denote the Lasso estimator defined based on the sample $\left(\xa,\ya\right)$ with the proper tuning parameter $\lambda=\sqrt{\frac{2.05 \log p}{n_1}}\sigma_0$,
\begin{equation}
\widehat{\beta}=\argmin_{\beta \in \R^{p}}\frac{\|\ya-\xa\beta\|_2^2}{2n_1}+\lambda \sum_{j=1}^{p} \frac{\|\xa_{\cdot j}\|_2}{\sqrt{n_1}} |\beta_j|.
\label{eq: Lasso estimator}
\end{equation}
We define the following estimator of $\xi^{\intercal} \beta$,
\begin{equation}
\bar{\mu}=\xi^{\intercal} \widehat{\beta}+\frac{1}{n_2} \xi^{\intercal} \left(\xb\right)^{\intercal}\left(\yb-\xb \widehat{\beta}\right).
\label{eq: known Sigma estimator}
\end{equation}
Based on the estimator, we construct the following confidence interval
\begin{equation}
\CIZk=\left[\bar{\mu}-1.01 \frac{\|\xi\|_2}{ \sqrt{n_2}}\zas \sigma_0,\bar{\mu}+1.01 \frac{\|\xi\|_2}{ \sqrt{n_2}}\zas \sigma_0 \right],
\label{eq: CI length Known Sigma}
\end{equation}
where $\alpha_0= \gamma_0 \alpha$ with $0<\gamma_0<1$. It will be shown in the supplement \citep{cai2015supplement} that the confidence interval proposed in \eqref{eq: CI length Known Sigma} has valid coverage and achieves the minimax length in \eqref{eq: minimax for known Sigma}.
\subsection{Confidence intervals for linear functionals in the dense loading regime}
\label{sec: dense functional known}
In marked contrast to the sparse loading regime, the prior knowledge of $\Omega={\rm I}$ and $\sigma=\sigma_0$ does not improve the minimax rate in the dense loading regime. That is,  Theorem \ref{thm: minimax for dense linear functional} remains true by replacing $\Theta\left(k\right)$ and $\Theta\left(k_1\right)$ with $\Theta\left(k,{\rm I},\sigma_0\right)$ and $\Theta\left(k_1,{\rm I},\sigma_0\right)$, respectively. However, the cost of adaptation changes when there is prior knowledge of $\Omega={\rm I}$ and $\sigma=\sigma_0$. 
The following theorem establishes the adaptivity lower bound in the dense loading regime.
\begin{Theorem}
\label{thm: adaptivity multi dense both known}
Suppose that $0 < \alpha <\frac{1}{2}$ and $k_1\le k \leq c\min\left\{p^{\gamma},\frac{n}{\log p}\right\}$ for some constants $c>0$ and $0\leq \gamma <\frac{1}{2}$, then, for some constant $c_1>0$,
\begin{equation}
\label{eq: adaptivity multi dense both known}
L_{\alpha}^{*}\left(\Theta\left({k_1},{\rm I},\sigma_0\right),\Theta\left({k},{\rm I},\sigma_0\right),\xi^{\intercal} \beta \right)\geq c_1 \|\xi\|_{\infty}\sigma_0\max\left\{\sqrt{k k_1} \sqrt{\frac{\log p}{n}},\min\left\{k \sqrt{\frac{\log p}{n}},\frac{\sqrt{k}}{n^{\frac{1}{4}}}\right\}\right\}.
\end{equation}
\end{Theorem}


The lower bound in \eqref{eq: adaptivity multi dense both known} is attainable. For reasons of space, the construction is omitted here.
Under the framework \eqref{eq: adaptive benchmark}, adaptive confidence intervals are still impossible, since for $k_1\ll k$,
$$ L_{\alpha}^{*}\left(\Theta\left({k_1},{\rm I},\sigma_0\right),\Theta\left({k},{\rm I},\sigma_0\right),\xi^{\intercal} \beta \right) \gg L_{\alpha}^{*}\left(\Theta\left({k_1},{\rm I},\sigma_0\right),\xi^{\intercal} \beta \right).$$
Compared with Theorem \ref{thm: adaptivity multi dense}, we observe that the cost of adaptation is reduced with the prior knowledge of $\Omega={\rm I}$  and $\sigma=\sigma_0$.

\section{Discussion}
\label{sec:discussion}

In the present paper we studied the minimaxity and adaptivity of confidence intervals for general linear functionals $\xi^{\intercal} \beta$ with a sparse or dense loading $\xi$ for the setting where $\Omega$ and $\sigma$ are unknown as well as the setting with the prior knowledge of $\Omega={\rm I}$ and $\sigma=\sigma_0$. In the more typical case in practice where  $\Omega$ and $\sigma$ are unknown, the adaptivity results are quite negative: With the exception of the ultra-sparse region for confidence intervals for  $\xi^{\intercal} \beta$ with a sparse loading $\xi$, it is necessary to know the true sparsity $k$ in order to have guaranteed coverage probability and rate-optimal expected length.  In contrast to estimation, knowledge of the sparsity $k$ is crucial to constructing honest confidence intervals.  In this sense,  the problem of constructing confidence intervals is much harder than the estimation problem.

The case of known $\Omega={\rm I}$ and $\sigma=\sigma_0$ is strikingly different. The minimax expected length in the sparse loading regime is of order $\frac{\|\xi\|_2}{\sqrt{n}}$ and in particular does not depend on $k$ and adaptivity can be achieved over the full range of sparsity $k\lesssim \frac{{n}}{\log p}$. So in this case, the knowledge of $\Omega$ and $\sigma$ is very useful.
On the other hand,  in the dense loading regime the information on $\Omega$ and $\sigma$ is of limited use.
In this case, the minimax rate and lack of adaptivity remain unchanged, compared with the unknown $\Omega$ and $\sigma$ case, although the cost of adaptation is reduced.


Regarding the construction of confidence intervals, there is a significant difference between the sparse and dense loading regimes. The de-biasing method is useful in the sparse loading regime since such a procedure reduces the bias but does not dramatically increase the variance. However, the de-biasing construction is not applicable to the dense loading regime since the cost of obtaining a near-unbiased estimator is to significantly increase the variance which would lead to an unnecessarily long confidence interval.  An interesting open problem is the construction of a confidence interval for $\xi^{\intercal} \beta$  achieving the minimax length where the sparsity $q$ of the loading $\xi$ is in the middle regime with $c p^{\gamma} \leq q \leq c p^{2\gamma+\varsigma}$ for some $0<\varsigma<1-2\gamma$. 

In addition to constructing confidence intervals for linear functionals, another interesting problem is constructing confidence balls for the whole vector $\beta$. Such has been considered in \citep{nickl2013confidence}, where the authors established the impossibility of adaptive confidence balls for sparse linear regression. 
These problems are connected, but each has its own special features and the behaviors of the problems are different from each other. The connections and differences in adaptivity among various forms of confidence sets have also been observed in nonparametric function estimation problems. See, for example,  \citep{cai2004adaptation}  for adaptive confidence intervals for linear functionals, \citep{hoffmann2011adaptive,cai2014adaptive} for adaptive confidence bands, and \citep{Cai06,robins2006adaptive} for adaptive confidence balls.



In the context of nonparametric function estimation, a general adaptation theory for confidence intervals for an arbitrary linear functional was developed in Cai and Low \citep{cai2004adaptation} over a collection of convex parameter spaces. It was shown that the key quantity that determines adaptivity is a geometric quantity called the between-class  modulus of continuity. The convexity assumption on the parameter space  in Cai and Low \citep{cai2004adaptation} is crucial for the adaptation theory.  In high-dimensional linear regression, the parameter space is highly non-convex. The adaptation theory developed  in \citep{cai2004adaptation} does not apply to the present setting of high-dimensional linear regression.  It would be of significant interest to develop a general adaptation theory for confidence intervals in such a non-convex setting.


\section{Proofs}
\label{sec:proof}

In this section, we prove two main results, Theorem \ref{thm: strong non-adaptivity multi sparse} and minimax upper bound of Theorem \ref{thm: minimax for sparse linear functional}. For reasons of space, the proofs of the other results are given in the supplement \citep{cai2015supplement}.

A key technical tool for the proof of the lower bound results is the following lemma  which establishes the adaptivity over two nested parameter spaces. Such a formulation has been considered in \citep{cai2004adaptation} in the context of adaptive confidence intervals over convex parameter spaces under the Gaussian sequence model. However, the parameter space $\Theta(k)$ considered in the high dimension setting is highly non-convex. The following lemma can be viewed as a generalization of \citep{cai2004adaptation} to the non-convex parameter space, where the lower bound argument requires testing for  composite hypotheses.


Suppose that we observe a random variable $Z$ which has a distribution $\mathbf{P}_{\theta}$ where the parameter $\theta$ belongs to the parameter space $\HH$. Let $\CIZg$ be the confidence interval for the linear functional $\T\left(\theta\right)$ with the guaranteed coverage $1-\alpha$ over the parameter space $\HH$. Let $\HH_{0}$ and $\HH_{1}$ be subsets of the parameter space $\HH$ where $\HH=\HH_{0} \cup \HH_{1}$. 
Let $\pi_{\HH_i}$ denote the prior distribution supported on the parameter space $\HH_{i}$ for $i=0,1$.  
Let $f_{\pi_{\HH_i}}\left(z\right)$ denote the density function of the marginal distribution of $Z$  with the prior $\pi_{\HH_{i}}$ on $\HH_{i}$ for $i=0,1$. More specifically,
$f_{\pi_{\HH_i}}\left(z\right)=\int f_{\theta}\left(z\right) \pi_{\HH_i}\left(\theta\right)d\theta.$ 

Denote by $\mathbb{P}_{\pi_{\HH_i}}$ the marginal distribution of $Z$ with the prior $\pi_{\HH_{i}}$ on $\HH_{i}$ for $i=0,1$.  For any function $g$, we write $\E_{{\pi_{\HH_0}}}\left(g\left(Z\right)\right)$ for the expectation of $g\left(Z\right)$ with respect to the marginal distribution of $Z$ with the prior $\pi_{\HH_{0}}$ on $\HH_{0}$. We define the $\chi^2$ distance between two density functions $f_{1}$ and $f_{0}$  by
\begin{equation}
\chi^2(f_{1},f_{0})=\int \frac{\left(f_1(z)-f_0(z)\right)^2}{f_0(z)} dz=\int \frac{f^2_{1}(z)}{f_{0}(z)} dz-1
\label{eq: chisq distance}
\end{equation}
and the total variation distance by
$\TV(f_{1},f_{0})=\int \left|f_1(z)-f_0(z)\right| dz.$
It is well known that 
\begin{equation}
\TV(f_{1},f_{0})\leq \sqrt{\chi^2(f_{1},f_{0})}.
\label{eq: relation between chisq and TV}
\end{equation}
\begin{Lemma}
\label{lem: adaptivity general}
Assume $\T\left(\theta\right)=\mu_0$ for $\theta \in \HH_0$ and $\T\left(\theta\right)=\mu_1$ for $\theta \in \HH_1$ and $\HH=\HH_{0} \cup \HH_{1}.$ For any $\CIZg \in \II_{\alpha}\left(\T,\HH\right)$,
\begin{equation}
L\left(\CIZg,\HH\right) \geq L\left(\CIZg,\HH_0\right) \geq |\mu_1-\mu_0| \left(1-2 \alpha-{\TV\left(f_{\pi_{\HH_1}},f_{\pi_{\HH_0}}\right)}\right)_{+}.
\label{eq: non-adaptive for the gen linear functional}
\end{equation}
\end{Lemma}
\subsection{Proof of Lemma \ref{lem: adaptivity general}}
\label{sec: proof lower adaptivity key lemma}
The supremum risk over $\HH_0$ is lower bounded by the Bayesian risk with the prior $\pi_{\HH_0}$ on $\HH_0$,
\begin{equation}
\sup_{\theta \in \HH_0}\E_{\theta} L\left(\CIZg\right) \geq \int_{\theta} \E_{\theta} L\left(\CIZg\right) \pi_{\HH_0}\left(\theta\right) d\theta=\E_{{\pi_{\HH_0}}} L\left(\CIZg\right).
\label{eq: max by bayesian gen}
\end{equation}
By the definition of $\CIZg \in \II_{\alpha}\left(\T,\HH\right),$ we have
\begin{equation}
\mathbb{P}_{{\pi_{\HH_i}}}\left(\mu_i \in \CIZg\right)=\int_{\theta} \mathbb{P}_{\theta}\left(\mu_i \in \CIZg \right) \pi_{\HH_i}\left(\theta\right) d\theta \geq 1-\alpha,
\label{eq: control of prob 1 gen}
\end{equation}
for $i=0,1$. By the following inequality
$$\left |\mathbb{P}_{{\pi_{\HH_1}}}\left(\mu_1 \in \CIZg\right)-\mathbb{P}_{{\pi_{\HH_0}}}\left(\mu_1 \in \CIZg\right)\right |
\leq \TV\left(f_{\pi_{\HH_1}},f_{\pi_{\HH_0}}\right),$$
then we have
$
\mathbb{P}_{{\pi_{\HH_0}}}\left(\mu_1 \in \CIZg\right) \geq 1-\alpha-\TV(f_{\pi_{\HH_1}},f_{\pi_{\HH_0}}).
$
This together with $(\ref{eq: control of prob 1 gen})$  yields
$\mathbb{P}_{{\pi_{\HH_0}}}\left(\mu_0, \mu_1 \in \CIZg\right) \geq 1-2\alpha-\TV(f_{\pi_{\HH_1}},f_{\pi_{\HH_0}}),$
which leads to
$\mathbb{P}_{{\pi_{\HH_0}}}\left(L\left(\CIZg\right)\geq \left|\mu_1-\mu_0\right|\right) \geq 1-2\alpha-\TV(f_{\pi_{\HH_1}},f_{\pi_{\HH_0}}).$
Hence,
$\E_{{\pi_{\HH_0}}} L (\CIZg) \geq \left(\mu_1-\mu_0\right) (1-2\alpha-\TV(f_{\pi_{\HH_1}},f_{\pi_{\HH_0}}))_{+}.$
The lower bound $(\ref{eq: non-adaptive for the gen linear functional})$ follows from inequality $(\ref{eq: max by bayesian gen}).$
\subsection{Proof of Theorem \ref{thm: strong non-adaptivity multi sparse}}
\label{sec: proof lower multi sparse non-parametric}

The lower bound in \eqref{eq: strong non-adaptivity multi sparse} is involved with a parametric term and a non-parametric term. The proof of the parametric lower bound is postponed to the supplement. In the following, we will prove the non-parametric lower bound
\begin{equation}
\inf_{\CIZl \in \SCIZl }\E_{\theta^*} L\left(\CIZl  \right)\geq c_1 \|\xi\|_2 k\frac{\log p}{n} \sigma,
\label{eq: the nonparametric goal}
\end{equation}
for some constant $c_1>0$.
Without loss of generality, we assume $\suppxi=\{1,\cdots,\|\xi\|_0\}$. 
We generate the orthogonal matrix $M \in \R^{\|\xi\|_0\times \|\xi\|_0}$ such that its first row is $\frac{1}{\|\xi\|_2} \xi_{\suppxi}$ and
define the orthogonal matrix $Q$ as $Q=\begin{pmatrix} M & \mathbf{0}\\ \mathbf{0}& \mathbf{\rm I} \end{pmatrix}$. 
We transform both the design matrix $X$ and the regression vector $\beta$ and view the linear model \eqref{eq: linear model} as
$y=V\psi+\epsilon,$
where $V=XQ^{\intercal}$ and $\psi=Q\beta$. The transformed coefficient vector $\psi^{*}=Q\beta^{*}=\begin{pmatrix} M \beta^{*}_{\suppxi}\\ \beta^{*}_{-\suppxi}\end{pmatrix}$ is of sparsity at most $\|\xi\|_0+k_1$.  The first coefficient $\psi_1$ of $\psi$ is $\frac{1}{\|\xi\|_2} \xi^{\intercal} \beta$.
The covariance matrix $\Psi$ of  $V_{1\cdot}$ is $Q \Sigma Q^{\intercal}$ and its corresponding precision matrix is  $\Gamma=Q\Omega Q^{\intercal}.$ 
To represent the transformed observed data and parameter, we abuse the notation slightly and also use $Z_i=\left(y_i,V_{i\cdot}\right)$ and $\theta^{*}=\left(\psi^{*},{\rm I}, \sigma \right)$. We define the parameter space $\GG\left(k\right)$ of $\left(\psi,\Gamma,\sigma\right)$ as
\begin{equation}
\GG\left(k\right)=\left\{\left(\psi,\Gamma,\sigma\right): \psi=Q\beta,\; \Gamma=Q \Omega Q^{\intercal}\; \text{for} \; \left(\beta,\Omega,\sigma\right)\in \Theta\left(k\right) \right\}.
\label{eq: parameter space of transformed para}
\end{equation}

For a given $Q$, there exists a bijective mapping between $\Theta\left(k\right)$ and $\GG\left(k\right)$. To show that $\left(\psi,\Gamma,\sigma\right) \in \GG\left(k\right)$, it is equivalent to show $\left(Q^{\intercal}\psi, Q^{\intercal} \Gamma Q,\sigma \right)\in \Theta\left(k\right)$.
Let $\SCIZtrans$ denote the set of confidence intervals for $\psi_1=\frac{1}{\|\xi\|_2} \xi^{\intercal} \beta$ with guaranteed coverage over $\GG\left(k\right)$. If $\CIZtrans\in \SCIZtrans$, then $\|\xi\|_2\CIZtrans \in \SCIZl$; If $\CIZl \in \SCIZl$, then $\frac{1}{\|\xi\|_2} \CIZl \in \SCIZtrans$.
Because of such one to one correspondence, we have
\begin{equation}
\inf_{\CIZl \in \SCIZl }\E_{\theta^*} L\left(\CIZl  \right) = \|\xi\|_2 \inf_{\CIZtrans \in \SCIZtrans }\E_{\theta^*} L\left(\CIZtrans \right).
\label{eq: equivalence of transformation}
\end{equation}
By \eqref{eq: the nonparametric goal} and \eqref{eq: equivalence of transformation}, we reduce the problem to 
\begin{equation}
\inf_{\CIZtrans \in \SCIZtrans }\E_{\theta^*} L\left(\CIZtrans \right) \geq c {k} \frac{\log p}{n}\sigma.
\label{eq: after transformation goal}
\end{equation}
Under the Gaussian random design model, $Z_{i}=\left(y_i,V_{i\cdot}\right)\in \R^{p+1}$
follows a joint Gaussian distribution with mean $0$. Let $\Sigma^z$ denotes the covariance matrix of $Z_{i}$. 
Decompose $\Sigma^{z}$ into blocks $\begin{pmatrix} \Sigma_{yy}^{z}& \left(\Sigma_{vy}^{z}\right)^{\intercal}\\ \Sigma_{vy}^{z}& \Sigma_{vv}^{z} \end{pmatrix},$ where $\Sigma_{yy}^{z}$, $\Sigma_{vv}^{z}$ and $\Sigma_{vy}^{z}$ denote the variance of $y$, the variance of $V$ and the covariance of $y$ and $V$, respectively. We define the function $h : \Sigma^z \rightarrow \left(\psi,\Gamma,\sigma\right)$  as  $h(\Sigma^z)=\left(\left(\Sigma_{vv}^{z}\right)^{-1}\Sigma_{vy}^{z}, \left(\Sigma_{vv}^{z}\right)^{-1},\Sigma_{yy}^{z}-\left(\Sigma_{vy}^{z}\right)^{\intercal}\left(\Sigma_{vv}^{z}\right)^{-1}\Sigma_{vy}^{z}\right)$. The function $h$ is bijective and  its inverse mapping $h^{-1}: \left(\psi,\Gamma,\sigma\right) \rightarrow \Sigma^z$ is $$h^{-1}\left(\left(\psi,\Gamma,\sigma\right)\right)=\begin{pmatrix} \psi^{\intercal}\Gamma^{-1}\psi+\sigma^2& \psi^{\intercal}\Gamma^{-1}\\ \Gamma^{-1}\psi& \Gamma^{-1}\end{pmatrix}.$$

The null space is taken as $\HH_0=\left\{\left(\psi^{*},{\rm I},\sigma \right)\right\}$ and  $\pi_{\HH_0}$ denotes the point mass prior at this point. The proof is divided into three steps:
\begin{enumerate}
\item Construct $\HH_1$ and show that $\HH_1\subset \GG\left(k\right)$; 
\item Control the distribution distance $\TV\left(f_{\pi_{\HH_1}},f_{\pi_{\HH_0}}\right)$;
\item Calculate the distance $\mu_1-\mu_0$ where $\mu_0=\psi_1^{*}$ and $\mu_1=\psi_1$ with $\left(\psi,\Gamma,\sigma\right)\in \HH_1$. We show that $\mu_1= \psi_1$ is a fixed constant for all  $\left(\psi,\Gamma,\sigma\right)\in \HH_1$ and  then apply Lemma \ref{lem: adaptivity general}.

 \end{enumerate}
\textbf{Step 1}. We construct the alternative hypothesis parameter space $\HH_1$.
Let $\Sigma^z_{0}$ denote the covariance matrix of $Z_i$ corresponding to $\left(\psi^{*},{\rm I},\sigma\right)\in \HH_0$.  Let $S_1={\rm supp}\left(\psi^{*}\right) \cup \left\{1\right\}$ and $S=S_1 \backslash\{1\}$. Let $k_{*}$ denote the size of $S$ and $p_1$ denote the size of $S_1^c$ and  we have $k_{*}\leq k_1+q$ and $p_1\geq p-k_{*}-1\geq c p$. 
Without loss of generality, let $S=\{2,\cdots,k_{*}+1\}$. We have the following expression for the covariance matrix of $Z_i$ under the null,
\begin{equation}
\Sigma_0^{z}=\left(
\begin{array}{c|c|c|c}
 \|\psi^*\|_2^2+\sigma^2&\psi_{1}^{*}&\left(\psi_{S}^{*}\right)^{\intercal}& \mathbf{0}_{1\times p_1}\\ \hline
\psi_{1}^{*}&1& \mathbf{0}_{1\times k_{*}}&\mathbf{0}_{1\times p_1}\\ \hline
 \psi_{S}^{*}&\mathbf{0}_{k_{*}\times 1}& \rm I_{k_{*}\times k_{*}}&\mathbf{0}_{k_{*}\times p_1}\\ \hline
 \mathbf{0}_{p_1\times 1}&\mathbf{0}_{p_1\times 1}&\mathbf{0}_{ p_1\times k_{*}}&{\rm I}_{p_1\times p_1}\\
\end{array}
\right),
\label{eq: null Z}
\end{equation}
To construct $\HH_1$, we define the following set, 
\begin{equation}
\ell\left(p_1,\frac{\zeta_0}{2} k,\rho\right)=\left\{\alphab :  \alphab \in \R^{p_1},\, \|\alphab\|_0=\frac{\zeta_0}{2} k, \,\alphab_{i} \in\{0,\rho\} \; \text{for} \;1\leq i\leq p_1\right\}.
\label{eq: prior of alternative}
\end{equation}
Define the parameter space $\FF$ for $\Sigma^{z}$ by
$ \FF=\left\{\Sigma^{z}_{\alphab}: {\alphab} \in \ell\left(p_1,\frac{\zeta_0}{2} k,\rho\right) \right\}$,
where
\begin{equation}
\Sigma_{\alphab}^{z}=\left(
\begin{array}{c|c|c|c}
 \|\psi^*\|_2^2+\sigma^2&\psi_1^{*}&\left(\psi_{S}^{*}\right)^{\intercal}& \rho_0 \alphab^{\intercal}\\ \hline
 \psi_1^{*}&1& \mathbf{0}_{1\times k_{*}}&\alphab^{\intercal}\\ \hline
 \psi_{S}^{*}&\mathbf{0}_{k_{*}\times 1}& \rm I_{k_{*}\times k_{*}}&\mathbf{0}_{k_{*}\times p_1}\\ \hline
\rho_0 \alphab&\alphab&\mathbf{0}_{ p_1\times k_{*}}&\rm I_{p_1\times p_1}\\
\end{array}
\right).
\label{eq: alter Z}
\end{equation}
Then we construct the alternative hypothesis space $\HH_1$ for $\left(\psi,\Gamma,\sigma\right)$, which is induced by the mapping $h$ and the parameter space $\FF$,
\begin{equation}
\HH_1=\left\{\left(\psi,\Gamma,\sigma\right) \; :\; \left(\psi,\Gamma,\sigma\right)=h\left(\Sigma^{z}\right) \quad \text{for} \; \Sigma^{z} \in \FF\right\}
\label{eq: def of alternative space}
\end{equation}
In the following, we show that $\HH_1 \subset \GG\left(k\right)$. 
It is necessary to identify  $\left(\psi,\Gamma,\sigma\right)=h\left(\Sigma^{z}\right)$ for $\Sigma^{z} \in \FF$ and show $\left(Q^{\intercal} \psi, Q^{\intercal}\Gamma Q,\sigma \right)\in \Theta\left(k\right)$.  Firstly, we identify the expression $\E\left(y_{i} \mid V_{i,\cdot}\right)$ under the alternative joint distribution \eqref{eq: alter Z}.
Assuming $y_{i}=V_{i1} \psi_1+ V_{i,S}\psi_{S}+V_{i,S_1^c}\psi_{S_1^c}+\epsilon'_i,$
we have
\begin{equation}
\psi_1=\frac{-\|\alphab\|_2^2 \rho_0+\psi_1^{*}}{1-\|\alphab\|_2^2}, \;\psi_{S}=\psi_{S}^*, \;  \psi_{S_1^c}=\left(\rho_0-\psi_1\right)\alphab,
\label{eq: solution for alter}
\end{equation}
and
\begin{equation}
{\rm Var}\left(\epsilon'_i\right)=\sigma^2-\frac{\|\alphab\|_2^2\left(\rho_0-\psi^{*}_1\right)^2}{1-\|\alphab\|_2^2} \leq \sigma^2 \leq \Mv.
\label{eq: variance alter}
\end{equation}
Based on $\eqref{eq: solution for alter}$, the sparsity of $\psi$ in the alternative hypothesis space is upper bounded by
$1+|{\rm supp}\left(\psi_{S}^*\right)|+|{\rm supp}\left(\alphab\right)|\leq \left(1-\frac{\zeta_0}{4}\right)k,$ and hence the sparsity of the corresponding $\beta=Q^{\intercal} \psi$ is controlled by
\begin{equation}
\|\beta\|_{0}\leq \left(1-\frac{\zeta_0}{4}\right)k+q\leq k.
\label{eq: original beta in space}
\end{equation}
Secondly, we show that $\Omega=Q^{\intercal}\Gamma Q$ satisfies the condition 
$\frac{1}{\Md} \leq \lambda_{\min}\left(\Omega\right)\leq \lambda_{\max}\left(\Omega\right) \leq \Md.$
The covariance matrix $\Psi$ of $V_{i,\cdot}$ in the alternative hypothesis parameter space is expressed as
\begin{equation}
\Psi
=\left(\begin{array}{c|c|c}
1& \mathbf{0}_{1\times k_*}&\mathbf{0}_{k_*\times p_1}\\ \hline
\mathbf{0}_{k_*\times 1}& \rm I_{k_*\times k_*}&\mathbf{0}_{k_*\times p_1}\\ \hline
\mathbf{0}_{p_1 \times 1}&\mathbf{0}_{p_1\times k_*}&\rm I_{p_1\times p_1}\\
\end{array}
\right)+\left(\begin{array}{c|c|c}
0& \mathbf{0}_{1\times k_*}&\alphab^{\intercal}\\ \hline
\mathbf{0}_{k_*\times 1}& \mathbf{0}_{k_*\times k_*}&\mathbf{0}_{k_*\times p_1}\\ \hline
\alphab&\mathbf{0}_{ p_1\times k_*}& \mathbf{0}_{p_1\times p_1}\\
\end{array}
\right).
\end{equation}
Since the second matrix on the above equation is of spectral norm $\|\alphab\|_2$, Weyl's inequality 
leads to
$\max\left\{ \left|\lambda_{\min}\left(\Psi\right)-1\right|,\left|\lambda_{\max}\left(\Psi\right)-1\right| \right\} \leq \|\alphab\|_2 .$
When $\|\alphab\|_2$ is chosen such that $\|\alphab\|_2 \leq \min\left\{1-\frac{1}{\Md},\Md-1\right\},$ then we have $\frac{1}{\Md} \leq \lambda_{\min}\left(\Psi\right)\leq \lambda_{\max}\left(\Psi\right) \leq \Md.$ Since $\Omega$ and $\Gamma=Q \Omega Q^{\intercal}$ have the same eigenvalues, we have $\frac{1}{\Md} \leq \lambda_{\min}\left(\Omega\right)\leq \lambda_{\max}\left(\Omega\right) \leq \Md$. Combined with \eqref{eq: variance alter} and \eqref{eq: original beta in space}, we show that $\HH_1\subset \GG\left(k\right)$. 

\medskip\noindent
\textbf{Step 2}. To control $\TV\left(f_{\pi_{\HH_1}},f_{\pi_{\HH_0}}\right)$, it is sufficient to control $\chi^2\left(f_{\pi_{\HH_1}},f_{\pi_{\HH_0}}\right)$ and apply \eqref{eq: relation between chisq and TV}. Let $\pi$ denote the uniform prior on $\alphab$ over $\ell\left(p_1,\frac{\zeta_0}{2}k,\rho\right)$. Note that this uniform prior $\pi$ induces a prior distribution $\pi_{\HH_1}$ over the parameter space $\HH_1$. 
Let $\E_{\alphab,\widetilde{\alphab}} $ denote the expectation with respect to the independent  random variables $\alphab,\widetilde{\alphab}$ with uniform prior $\pi$ over the parameter space $\ell\left(p_1,\frac{\zeta_0}{2} k,\rho\right)$.
The following lemma controls the $\chi^2$ distance between the null and the mixture over the alternative distribution.
\begin{Lemma}
\label{lem: workhorse1}
Let $f_1=\left(\sigma^2+\left(\psi_1^*\right)^2-\rho_0\psi_1^*\right)$. Then
\begin{equation}
\label{eq: control of chisq for beta1}
\chi^2\left(f_{\pi_{\HH_1}},f_{\pi_{\HH_0}}\right)+1= \E_{\alphab,\widetilde{\alphab}} \left(1-\frac{1}{\sigma^2} \left(\rho_0\left(\rho_0-\psi_1^{*}\right)+f_1 \right) \alphab^{\intercal}\widetilde{\alphab}\right)^{-n}.
\end{equation}
\end{Lemma}
The following lemma is useful in controlling the right hand side of \eqref{eq: control of chisq for beta1}.
\begin{Lemma}
\label{lem: bound the MGF of Hyper}
Let $J$ be a ${\rm Hypergeometric}\left(p,k,k\right)$ variable with
$\mathbb{P}\left(J=j\right)=\frac{\binom{k}{j} \binom{p-k}{k-j}}{\binom{p}{k}},$
then
\begin{equation}
\label{eq: upper bound for MGF of Hyper}
\E\exp\left(t J\right)\leq e^{\frac{k^2}{p-k}}\left(1-\frac{k}{p}+\frac{k}{p}\exp\left(t\right)\right)^{k}.
\end{equation}
\end{Lemma}
Taking $\rho_0=\psi_1^*+\sigma$, we have
$\frac{1}{\sigma^2} \left(\rho_0\left(\rho_0-\psi_1^{*}\right)+f_1 \right)=2$ and by Lemma \ref{lem: workhorse1},
\begin{equation*}
\chi^2\left(f_{\pi_{\HH_1}},f_{\pi_{\HH_0}}\right)+1=\E_{\alphab,\widetilde{\alphab}} \left(1-2\alphab^{\intercal}\widetilde{\alphab}\right)^{-n}.
\end{equation*}
By the inequality $\frac{1}{1-x}\leq \exp(2x)$ for $x\in \left[0,\frac{\log 2}{2}\right]$,
if $\alphab^{\intercal} \widetilde{\alphab} \leq \frac{\zeta_0}{2} k\rho^2 <\frac{\log 2}{4},$ then $\left(1-2\alphab^{\intercal} \widetilde{\alphab}\right)^{-n}\leq  \exp\left( 4 n \alphab^{\intercal} \widetilde{\alphab}\right)$. 
By Lemma \ref{lem: bound the MGF of Hyper}, we further have
\begin{equation*}
\begin{aligned}
&\E_{\alphab,\widetilde{\alphab}}  \exp\left( 4 n \alphab^{\intercal} \widetilde{\alphab}\right)= \E \exp\left(4J n\rho^2\right) \leq e^{\frac{\zeta_0^2k^2}{4p_1-2\zeta_0 k}}\left(1-\frac{\zeta_0 k}{2p_1}+\frac{\zeta_0 k}{2p_1} \exp\left(4 n\rho^2\right)\right)^{\frac{\zeta_0}{2} k}\\
 &\leq e^{\frac{\zeta_0^2 k^2}{4p_1-2\zeta_0 k}}\left(1-\frac{\zeta_0 k}{2 p_1}+\frac{\zeta_0 k}{2p_1} \sqrt{\frac{4p_1}{\zeta_0^2 k^2}}\right)^{\frac{\zeta_0}{2} k}
 \leq e^{\frac{c^2 \zeta_0^2 p^{2\gamma}}{4p_1-2 c \zeta_0 p^{\gamma}}}\left(1+\frac{1}{\sqrt{p_1}}\right)^{\frac{c\zeta_0}{2} p^{\gamma}},\\
 \end{aligned}
 \end{equation*}
where the second inequality follows by plugging in $\rho= \sqrt{\frac{\log \frac{4p_1}{\zeta_0^2 k^2}}{8n}}$ and the last inequality follows by $k\leq c p^{\gamma}$. If $k\leq c\left\{\frac{n}{\log p}, p^{\gamma}\right\}$, where $0\leq \gamma <\frac{1}{2}$ and $c$ is a sufficient small positive constant, then
$k \rho^2<\min\left\{\frac{\log 2}{2\zeta_0},\left(1-\frac{1}{\Md}\right)^2,1\right\}$
and hence
\begin{equation}
\chi^2\left(f_{\pi_{\HH_1}},f_{\pi_{\HH_0}}\right)\leq \left(\frac{1}{2}-\alpha\right)^2 \quad \text{and} \quad \TV\left(f_{\pi_{\HH_1}},f_{\pi_{\HH_0}}\right)\leq \frac{1}{2}-\alpha.
\label{eq: control of chisq in sparse beta}
\end{equation}

\medskip\noindent
\textbf{Step 3}.   We calculate the distance between $\mu_1$ and $\mu_0$.
Under $\HH_0$, $\mu_0=\psi_1^*.$
Under $\HH_1$, $\mu_1=\psi_{1}=\frac{-\|\alphab\|_2^2 \rho_0+ \psi_1^{*}}{1-\|\alphab\|_2^2}.$
For $\alphab\in \ell\left(p_1,\frac{\zeta_0}{2} k,\rho\right)$,
$\|\alphab\|_2^2=\frac{\zeta_0}{2} k\rho^2$ and $\mu_1=\psi_{1}=\frac{-\frac{\zeta_0}{2} k\rho^2\left(\psi_1^{*}+\sigma\right)+\psi_1^{*} }{1-\frac{\zeta_0}{2} k\rho^2 }.$
Since $\rho$ is selected as fixed, $\mu_1=\psi_1$ is a fixed constant for $\left(\psi,\Omega,\sigma\right)\in \HH_1.$
Note that $\mu_1-\mu_0=\frac{\|\alphab\|_2^2\left(\psi_1^{*}- \rho_0\right)}{1-\|\alphab\|_2^2}=\frac{-\sigma \|\alphab\|_2^2}{1-\|\alphab\|_2^2},$
and it follows that
$
|\mu_1-\mu_0|= \sigma \frac{\|\alphab\|_2^2}{1- \|\alphab\|_2^2}\geq  c k \frac{\log \frac{4 p_1}{\zeta_0^2 k^2}}{n} \sigma.
$
Combined with \eqref{eq: relation between chisq and TV} and \eqref{eq: control of chisq in sparse beta}, Lemma \ref{lem: adaptivity general} leads to \eqref{eq: after transformation goal}.
By \eqref{eq: equivalence of transformation}, we establish \eqref{eq: strong non-adaptivity multi sparse}.
\subsection{Proof of upper bound in Theorem \ref{thm: minimax for sparse linear functional}}
The following propostion establishes the coverage property and the expected length of the constructed confidence interval constructed in \eqref{eq: constructed CI for multi sparse}. Such a confidence interval achieves the minimax length in \eqref{eq: minimax for multi sparse}.
\begin{Proposition}
Suppose that
$k \leq c_{*} \frac{{n}}{\log p},$ where $c_{*}$ is a small positive constant, then
\begin{equation}
\liminfnp \inf_{\theta \in \Theta\left(k\right)}\PP\left(\xi^{\intercal} \beta \in \CIZs\right)
                \geq 1-\alpha,
\label{eq: covering property for multi sparse}
\end{equation}
and
\begin{equation}
L\left(\CIZs,\Theta\left(k\right) \right) \leq  C\|\xi\|_2 \left(k \frac{\log p}{n}+\frac{1}{\sqrt{n}}\right),
\label{eq: control of length for multi sparse}
\end{equation}
for some constant $C>0$.
\label{prop: covering thm multi sparse}
\end{Proposition}
In the following, we are going to prove Proposition \ref{prop: covering thm multi sparse}.
By normalizing the columns of $X$ and the true sparse vector $\beta$, the linear regression model can be expressed as
\begin{equation}
y=W d+\epsilon, \quad \text{with } \; W=X D, \; d=D^{-1}\beta \;\text{and} \; \epsilon \sim N(0,\sigma^2 {\rm I}),
\label{eq: normalized model}
\end{equation}
where
\begin{equation}
D={\rm diag}\left(\frac{\sqrt{n}}{\|\Xj\|_2}\right)_{j\in [p]}
\label{eq: def of D}
\end{equation}
 denotes the $p\times p$ diagonal matrix with $\left(j,j\right)$ entry to be $\frac{\sqrt{n}}{\|\Xj\|_2}$.
Take $\delta_0=1.0048$ and $\eta_0=0.01$, and we have $\lambda_0=\left(1+\eta_0\right)\sqrt{\frac{2\delta_0\log p}{n}}.$ Take $\epsilon_0=\frac{2.01}{\eta_0}+1=202$, $\nu_0=0.01$, $C_1=2.25$, $c_0=\frac{1}{6}$ and $C_0=3$. Rather than use the constants directly in the following discussion, we use $\delta_0,\eta_0,\epsilon_0,\nu_0, C_1,C_0$  and $c_0$ to represent the above fixed constants in the following discussion.  We also assume that
$\frac{\log p}{n} \leq \frac{1}{25}$ and $\delta_0\log p>2.$
Define the $l_1$ cone invertibility factor ($CIF_1$) as follows,
\begin{equation}
CIF_1\left(\alpha_0,K,W\right)=\inf\left\{\frac{|K| \|\frac{W^{\intercal} W}{n} u\|_{\infty}}{\|u_{K}\|_1}:\|u_{K^c}\|_1\leq \alpha_0 \|u_{K}\|_1, u\neq 0 \right\},
\end{equation}
where $K$ is an index set. Define
$\sigma^{ora}=\frac{1}{\sqrt{n}}\|y-X\beta\|_2=\frac{1}{\sqrt{n}}\|y-W d\|_2,$
\begin{equation}
T=\{k: |d_{k}|\geq \lambda_0 \sigma^{ora} \}, \; \tau=\left(1+\epsilon_0\right)\lambda_0\max \left\{\frac{4}{\sigma^{ora}}\|d_{T^{c}}\|_1,\frac{8 \lambda_0 |T|}{CIF_{1}\left(2 \epsilon_0+1,T,W\right)}\right\}.
\label{eq: def of T and tau}
\end{equation}
To facilitate the proof, we define the following events for the random design $X$ and the error $\epsilon$,
\allowdisplaybreaks
\begin{align*}
G_1&=\left\{ \frac{2}{5} \frac{1}{\sqrt{\Md}}<\frac{\|\Xj\|_2}{\sqrt{n}}<\frac{7}{5}\sqrt{\Md} \,\text{for}\,1\leq j\leq p\right\},\\
G_2&=\left\{ \left|\frac{\left(\sigma^{ora}\right)^2}{\sigma^2}-1\right|\leq 2\sqrt{\frac{\log p}{n}}+2\frac{\log p}{n}\right\},\\
G_3&=\left\{ \max\left\{\left|\frac{\xi^{\intercal} \widehat{\Sigma}\xi}{\xi^{\intercal} \Sigma \xi}-1\right|, \left|\frac{{\uu}^{\intercal} \widehat{\Sigma}{\uu}}{\xi^{\intercal} \Omega\xi}-1\right|\right\}\leq 2\sqrt{\frac{\log p}{n}}+2\frac{\log p}{n}\right\}, \; \text{where} \; \uu=\Omega \xi,\\
G_4&=\left\{\kappa(X,k,\alpha)\geq \frac{1}{4\sqrt{\lambda_{\max}\left(\Omega\right)}}- \frac{9}{\sqrt{\lambda_{\min}\left(\Omega\right)}}\left(1+\alpha\right) \sqrt{k \frac{\log p}{n}}\right\},\\
G_5&=\left\{\frac{\|W^{T} \epsilon\|_{\infty}}{n} \leq \sigma\sqrt{\frac{2 \delta_0 \log p}{n}}\right\},\\
S_1&=\left\{\frac{\|W^{T} \epsilon\|_{\infty}}{n} \leq \sigma^{ora} \lambda_0 \frac{\epsilon_0-1}{\epsilon_0+1}\left(1-\tau\right)\right\},\\
S_2&=\left\{\left(1-\nu_0\right)\hat{\sigma} \leq \sigma \leq (1+\nu_0) \hat{\sigma}\right\},\\
B_1&=\left\{\|\xi^{\intercal} \Omega\hat{\Sigma} - \xi^{\intercal}\|_{\infty}\leq \lambda_n\right\}, \; \text{ where } \;\lambda_n=4 C_0 \Md^2 \|\xi\|_2 \sqrt{\frac{\log p}{n}}.
\end{align*}

Define
$G=\cap_{i=1}^{5} G_i$ and $S=\cap_{i=1}^{2} S_i.$
The following lemmas  control the probability of events $G$, $S$ and $B_1$. The detailed proofs of Lemma \ref{lem: high prob}, \ref{lem: upper bound for bias} and \ref{lem: control the length} are in the supplement. 
\begin{Lemma}
\begin{equation}
\PP\left(G\right)\geq 1- \frac{6}{p}- 2p^{1-C_1}-\frac{1}{2\sqrt{\pi \delta_0 \log p}} p^{1-\delta_0}- c'\exp\left(-cn\right),
\label{eq: high prob}
\end{equation}
and
\begin{equation}
\PP\left(B_1\right) \geq 1-2p^{1-c_0 C^2_0},
\label{eq: high prob tuning}
\end{equation}
where $c$ and $c'$ are universal positive constants.
If $k \leq c_{*} \frac{n}{\log p}$, then
\begin{equation}
\PP\left(G \cap S\right)\geq \PP\left(G\right)-2 \exp\left(-\left(\frac{g_0+1-\sqrt{2g_0+1}}{2}\right)n\right)- c'' \frac{1}{\sqrt{\log p}} p^{1-\delta_0}, \label{eq: high prob large n}
\end{equation}
where $c_{*}$ and $c''$ are universal positive constants and $g_0=\frac{\nu_0}{2+3 \nu_0}$.
\label{lem: high prob}
\end{Lemma}
The following lemma establishes a data-dependent upper bound for the term $\|\widehat{\beta}-\beta\|_1$. 
\begin{Lemma}
On the event $G \cap S$,
\begin{equation}
\|\widehat{\beta}-\beta\|_1\leq \left(2+2\epsilon_0\right) \frac{\sqrt{n}}{\min\|\Xj\|_2} l\left(Z,k\right),
\label{eq: upper bound for bias}
\end{equation}
where
\begin{equation}
l\left(Z,k\right)=
\max\left \{ k\lambda_0\sigma^{ora},\frac{\left(2+2\epsilon_0\right) \max \|\Xj\|_2^2\left(\sigma \sqrt{\frac{2\delta_0\log p}{n}}+\lambda_0 \hat{\sigma}\right)k}{{n \kappa^2\left(X,k,\left(1+2\epsilon_0\right)\left(\frac{\max \|\Xj\|_2}{\min \|\Xj\|_2}\right)\right) }} \right \}.
\label{eq: def of upper bound}
\end{equation}
\label{lem: upper bound for bias}
\end{Lemma}
The following lemma controls the radius of the confidence interval.
\begin{Lemma}
\label{lem: control the length}
On the event $G\cap S \cap B_1$, there exists $p_0$ such that if $p\geq p_0$,
\begin{equation}
\rho_1\left(k\right) \leq C \|\xi\|_2 \left(\frac{1}{\sqrt{n}}+\frac{k \log p}{n}\right){\sigma}
\leq \|\xi\|_2 \log p \left(\frac{1}{\sqrt{n}}+\frac{k \log p}{n}\right)\hat{\sigma},\\
\label{eq: control the length3}
\end{equation}
and
\begin{equation}
\rho_2\left(k\right)\leq C k \sqrt{\frac{\log p}{n}} \sigma \leq \log p\left(k \sqrt{\frac{\log p}{n}}\hat{\sigma}\right).
\label{eq: control the length2}
\end{equation}
\end{Lemma}

In the following, we establish the coverage property of the proposed confidence interval.
By the definition of $\widetilde{\mu}$ in \eqref{eq: def of debiased estimator of linear functional}, we have
\begin{equation}
\widetilde{\mu}-\xi^{\intercal}\beta=\frac{1}{{n}}\widehat{\uu}^{\intercal} X^{\intercal} \epsilon+\left(\xi^{\intercal}-\widehat{\uu}^{\intercal}\widehat{\Sigma}\right)\left(\widehat{\beta}-\beta\right).
\label{eq: foundation of confidence interval}
\end{equation}
We now construct a confidence interval for the variance term $\frac{1}{{n}} \widehat{\uu}^{\intercal} X^{\intercal} \epsilon$ by normal distribution and a high probability upper bound for the bias term $ \left(\xi^{\intercal}-\widehat{\uu}^{\intercal}\widehat{\Sigma}\right)\left(\widehat{\beta}-\beta\right)$.
Since $\epsilon$ is independent of $X$ and $\widehat{\uu}$ and $\widehat{\Sigma}$ is a function of $X$, we have
$\frac{1}{{n}} \widehat{\uu}^{\intercal} X^{\intercal} \epsilon \mid X  \sim N\left(0,\sigma^2 \frac{\widehat{\uu}^{\intercal}\widehat{\Sigma}\widehat{\uu}}{n}\right),$
and
$$\mathbb{P}_{\epsilon\mid X } \left(\frac{1}{{n}} \widehat{\uu}^{\intercal} X^{T} \epsilon \in \left(-\sqrt{\frac{\widehat{\uu}^{\intercal}\widehat{\Sigma}\widehat{\uu}}{n}}\sigma \za,\sqrt{\frac{\widehat{\uu}^{\intercal}\widehat{\Sigma}\widehat{\uu}}{n}}\sigma \za \right) \middle| X \right)=1-\alpha.$$
By \eqref{eq: foundation of confidence interval}, we have
$\mathbb{P}_{\epsilon\mid X} \left(\xi^{\intercal}\beta \in {\rm CI}_0\left(Z,k\right)\middle| X  \right)=1-\alpha,$
where
\begin{equation*}
\begin{aligned}
{\rm CI}_0(Z,k) = \bigg[\widetilde{\mu}-\left(\xi^{\intercal}-\widehat{\uu}^{\intercal}\widehat{\Sigma}\right) &\left(\widehat{\beta}-\beta\right) -\sqrt{\frac{\widehat{\uu}^{\intercal}\widehat{\Sigma}\widehat{\uu}}{n}}\sigma \za, \\
& \widetilde{\mu} -\left(\xi^{\intercal}-\widehat{\uu}^{\intercal}\widehat{\Sigma}\right)\left(\widehat{\beta}-\beta\right)+ \sqrt{\frac{\widehat{\uu}^{\intercal}\widehat{\Sigma}\widehat{\uu}}{n}}\sigma \za \bigg].
\end{aligned}
\end{equation*}
Integrating with respect to $X$, we have
\begin{equation}
\PP \left(\xi^{\intercal} \beta \in {\rm CI}_0\left(Z,k\right)\right)=\int \mathbb{P}_{\epsilon\mid x} \left(\xi^{\intercal} \beta \in {\rm CI}_0\left(Z,k\right)\middle| x  \right) f(x) dx=1-\alpha.
\label{eq: oracle coverage}
\end{equation}
Since $\left|\left(\xi^{\intercal}-\widehat{\uu}^{\intercal}\widehat{\Sigma}\right)\left(\widehat{\beta}-\beta\right)\right|
\leq \|\xi^{\intercal}-\widehat{\uu}^{\intercal}\widehat{\Sigma}\|_{\infty}\|\widehat{\beta}-\beta\|_1,$
on the event $S \cap G$, Lemma \ref{lem: upper bound for bias} and the constraint in \eqref{eq: DZ for linear functional} lead to
\begin{equation}
\|\xi^{\intercal}-\widehat{\uu}^{\intercal}\widehat{\Sigma}\|_{\infty}\|\widehat{\beta}-\beta\|_1
\leq \lambda_n \left(2+2\epsilon_0\right) \frac{\sqrt{n}}{\min \|\Xj\|_2}l\left(Z,k\right),
\label{eq: upper bound for risk multi sparse}
\end{equation}
where $l\left(Z,k\right)$ is defined in \eqref{eq: def of upper bound}.
On the event $G \cap S$, we also have
$\sigma\leq \left(1+\nu_0\right) \hat{\sigma}$ and
$\sigma^{ora}\leq \left(1+\nu_0\right) \sqrt{1+2\sqrt{\frac{\log p}{n}}+2\frac{\log p}{n}} \hat{\sigma}.$
We define the following confidence interval to facilitate the discussion,
${\rm CI}_1\left(Z,k\right)
=\left(\widetilde{\mu}-l_k, \widetilde{\mu}+l_k\right),$
where $l_{k}=\left(1+\nu_0\right) \sqrt{\frac{\widehat{\uu}^{\intercal}\widehat{\Sigma}\widehat{\uu}}{n}}\za\hat{\sigma}+C_1\left(X,k\right)\|\xi\|_2 k \frac{\log p}{n} \hat{\sigma}.$
On the event $G \cap S$, we have
\begin{equation}
{\rm CI}_0\left(Z,k\right)\subset {\rm CI}_1\left(Z,k\right).
\label{eq: subset relation1}
\end{equation}
On the event $S_2$, if $p\geq \exp\left(2 \Mv\right)$, then $\hat{\sigma}\leq \frac{1}{1-\nu_0} \sigma\leq \frac{1}{1-\nu_0} \Mv<\log p$. Hence, the event $A$ holds and $\CIZs= \left[\widetilde{\mu}-\rho_1(k), \quad \widetilde{\mu}+\rho_1(k) \right]$.
By 
Lemma \ref{lem: control the length}, on the event $G \cap S \cap B_1$, if $p\geq \max\left\{p_0,\exp\left(2 \Mv\right)\right\}$, we have $\rho_1\left(k\right)=l_k,$
and hence 
\begin{equation}
{\rm CI}_1\left(Z,k\right)=\CIZs.
\label{eq: subset relation}
\end{equation}
We have the following bound on the coverage probability,
\begin{equation*}
\begin{aligned}
& \PP\left(\left\{\xi^{\intercal} \beta \in \CIZs\right\}\right)\geq \PP \left(\left\{\xi^{\intercal} \beta \in {\rm CI}_{0}\left(Z,k\right)\right\}\cap S \cap G\cap B_1\right)\\
\geq &\PP \left(\left\{\xi^{\intercal} \beta \in {\rm CI}_{0}\left(Z,k\right)\right\}\right)-\PP\left(\left(S \cap G\cap B_1\right)^c\right)
=1-\alpha-\PP\left(\left(S \cap G\cap B_1\right)^c\right)\\
=&\PP\left(S \cap G \cap B_1\right)-\alpha,
\end{aligned}
\end{equation*}
where the first inequality follows from \eqref{eq: subset relation1} and \eqref{eq: subset relation} and the first equality follows from \eqref{eq: oracle coverage}.
Combined with Lemma \ref{lem: high prob}, we establish \eqref{eq: covering property for multi sparse}.
We control the expected length as follows,
\begin{equation}
\begin{aligned}
&\E_{\theta} L\left(\CIZs \right) = \E_{\theta} L\left(\CIZs \right) \mathbf{1}_{A}\\
=&\E_{\theta} L\left(\CIZs \right) \mathbf{1}_{A\cap \left(S\cap G \cap B_1\right)}+\E_{\theta} L\left(\CIZs\right) \mathbf{1}_{A\cap \left(S\cap G \cap B_1\right)^c}\\
\leq & C\|\xi\|_2\left(k\frac{\log p}{n}+\frac{1}{\sqrt{n}}\right)\sigma+\|\xi\|_2 \left(\log p\right)^2 \left(\frac{1}{\sqrt{n}}+\frac{k \log p}{n}\right)\PP\left(\left(S\cap G \cap B_1\right)^c\right)\\
\leq &  C\|\xi\|_2\left(k\frac{\log p}{n}+\frac{1}{\sqrt{n}}\right) \left(\sigma+C \left(p^{1-\min\left\{\delta_0,C_1,c_0C_0^2\right\}}+c'\exp\left(-cn\right)\right)\left(\log p\right)^2\right),\\
\end{aligned}
\end{equation}
where the first inequality follows from \eqref{eq: control the length3} and second inequality follows from Lemma \ref{lem: high prob}. If $\frac{\log p}{n}\leq c$, then
$\left(p^{1-\min\left\{\delta_0,C_1,c_0C_0^2\right\}}+c'\exp\left(-cn\right)\right)\left(\log p\right)^2 \rightarrow 0,$
and hence
$
\E_{\theta} L\left(\CIZs \right)
\leq  C\|\xi\|_2\left(k\frac{\log p}{n}+\frac{1}{\sqrt{n}}\right)\Mv.
$